\documentclass[a4paper,12pt]{article}

\usepackage{listings}
\usepackage{amssymb,amsmath,amsfonts}
\usepackage[T1]{fontenc}
\usepackage{color}
\usepackage{amsthm}
\theoremstyle{remark}
\newtheorem{rem}{Remark}

\usepackage[width=16cm,height=26cm]{geometry}

\makeatletter
\newif\if@restonecol
\makeatother

\usepackage[algo2e,ruled,lined,titlenumbered,commentsnumbered]{algorithm2e}
\usepackage{caption}
\usepackage{subcaption}
\usepackage{mathtools}

\newcommand{\G}{\ensuremath{^G}}
\newcommand{\GT}{\ensuremath{^{G^T}}}
\newcommand{\Ls}{\ensuremath{^{s,L}}}
\newcommand{\Gs}{\ensuremath{^{s,G}}}

\newcommand{\Gz}{\ensuremath{^{0,G}}}

\newcommand{\LsT}{\ensuremath{^{s,L^T}}}

\newcommand{\s}{\ensuremath{^{s}}}

\newcommand{\z}{\ensuremath{^{0}}}
\newcommand{\sT}{\ensuremath{^{s^T}}}

\newcommand{\bK}{\ensuremath{\mathbf{K}}}

\newcommand{\bA}{\ensuremath{\mathbf{A}}}

\newcommand{\bT}{\ensuremath{\mathbf{T}}}

\newcommand{\bV}{\ensuremath{\mathbf{V}}}
\newcommand{\bW}{\ensuremath{\mathbf{W}}}

\newcommand{\bu}{\ensuremath{\mathbf{u}}}

\newcommand{\bp}{\ensuremath{\mathbf{p}}}

\newcommand{\br}{\ensuremath{\mathbf{r}}}

\newcommand{\foint}{\ensuremath{\mathbf{f}_{int}}}
\newcommand{\fointi}{\ensuremath{\mathbf{f}_{int,i}}}
\newcommand{\fointb}{\ensuremath{\mathbf{f}_{int,b}}}
\newcommand{\foext}{\ensuremath{\mathbf{f}_{ext}}}
\newcommand{\foexti}{\ensuremath{\mathbf{f}_{ext,i}}}
\newcommand{\foextb}{\ensuremath{\mathbf{f}_{ext,b}}}

\newcommand{\lam}{\ensuremath{\boldsymbol{\lambda}}}
\newcommand{\brho}{\ensuremath{\boldsymbol{\rho}}}

\newcommand{\brx}{\ensuremath{\mathbf{r}}}

\title{On the Implementation in Abaqus of the Global-Local Iterative Coupling and Acceleration Techniques}

\author{Omar Bettinotti$^1$, St\'ephane Guinard$^2$, Eric V\'eron$^1$ and Pierre Gosselet$^3$\\
{$^1$Dassault Syst\`emes Simulia Corp., Johnston, RI (USA)}\\
{$^2$IRT Saint Exupéry, Toulouse, France}\\
{$^3$Univ. Lille, CNRS, Centrale Lille, UMR 9013 - LaMcube, France}}

\begin{document}
\maketitle
            
\begin{abstract}
This paper presents results and convergence study of the Global-Local Iterative Coupling through the implementation in the commercial software Abaqus making use of the co-simulation engine. A hierarchical modeling and simulation approach is often required to alleviate modeling burdens. Particular focus has been devoted here on convergence acceleration and performance optimization. Two applications in statics with nonlinear material behavior and geometrically nonlinear formulation are considered here: first a holed curved plate under traction with elastic-plastic material, then a pre-stressed bolted joint connecting two plates between each other and subjected to traction load. Three different convergence acceleration techniques are compared in terms of convergence performance and accuracy. An inexact solver strategy is proposed to improve computing time performance. The results show promising results for the coupling technology and constitute a step forward in the availability of non-intrusive multi-scale modeling capabilities for complex structures and assemblies.
\end{abstract}



\section{Introduction}
\label{sec:intro}
Due to the complexity of aeronautical structure assemblies, engineers are used to building models across multiple levels of abstraction. At the highest level, a comprehensive model of the full aircraft is built from mid-surfacing and de-featuring the exact shapes of the assembly, with imposing a coarse finite elements mesh size, and simplifying assembly connections and interactions.
Material non-linearities are usually not considered at this level. The lower levels of abstractions are locally and hierarchically built, considering the exact geometry at the lowest level where a more refined finite elements mesh and material non-linearities are usually considered. As the lowest levels are local, their simulations are based on internal forces or boundary conditions coming from their highest level.

A main challenge today for aeronautic industry is to foster more agile connections between the different abstraction levels: by setting up accelerated paths from top to bottom levels, upwards and backwards (accelerated design); and also by refurbishing data flows between tests and simulations at different levels in the ``test pyramid'' (simulation aided certification).

As a preliminary step towards such high level objectives, the sub-modeling technique is often employed. It consists in applying the results of the global coarse simulation as boundary conditions in the local more refined simulation.
As the sub-modeling technique is sequential (global-to-local), an error is introduced and quantified as the difference between the integrals of the local reaction forces and global interface forces, where the local model cuts the global one.
This error is due to different modeling assumptions in the global and local models (constitutive law, connection modeling, mesh refinement\dots).
A more accurate coupling technique would require the global and local models to exchange data in both directions, paying the price of a higher computational cost, as simulations run concurrently.

Multiple strategies exist to couple global and local simulations together.
The most straightforward strategy consists in directly installing the local finite elements model into the global model with bonding interactions, possibly from solid to shell,  and then running a full monolithic simulation. An alternative strategy consists in using domain decomposition techniques, pioneered in \cite{lions_1990, glowinski_letallec_1990} as Schwarz alternating method, reformulated exchanging Neumann conditions in \cite{mandel_1993} or imposing the interface constraints by Lagrange multipliers in \cite{farhat_roux_1991, farhat_crivelli_roux_1994}.
A more comprehensive summary of domain decomposition techniques is presented in \cite{gosselet_rey_2006}.
A chief advantage of domain decomposition techniques is introduced in \cite{gravouil_combescure_2001, combescure_gravouil_2002, gravouil_combescure_2003, faucher_combescure_2003} for dynamic applications, as different time integration schemes or time increments can be applied to each sub-domain (sub-cycling), without compromising stability.
An interesting generalization to a multi-physic application is presented in  \cite{confalonieri_corigliano_gornati_dossi_2012}.

However, the chief issue with the aforementioned structural multi-scale techniques, the monolithic and the domain decomposition ones, is that the models abstractions need to be generated and adjusted for each simulation, aggravating the proliferation of models and increasing modeling complexity.
A less intrusive approach from the (non-negligible) modeling point of view is to keep the global coarse models unchanged, as it is the case of the sequential sub-modeling technique, and concurrently iterate between the global and local simulations. This is what we will refer to as the Global-Local Iterative Coupling (GLIC) in the following.
After the pioneering works of \cite{whitcomb_1991}, the formulation of the GLIC applied to nonlinear statics was revisited in \cite{gendre_allix_gosselet_comte_2009, gendre_allix_gosselet_2011}, reviewed in \cite{plews_duarte_eason_2012}, studied in the context of shell-to-solid coupling in \cite{guguin_gosselet_allix_guinard_2014, Guguin2016} and \cite{Akterskaia_2018} and of uncertainty quantification in \cite{nouy_pled_2018}.
In similar manner to the sub-cycling in domain decomposition, the GLIC was also extended to explicit dynamics in \cite{bettinotti_allix_malherbe_2013, bettinotti_allix_perego_oancea_malherbe_2014} and applied to a high-velocity impact causing composite delamination in \cite{bettinotti_2017}.
More recent works have proposed the application of the Generalized Finite Element Method to the local model in \cite{li_2021, li_2022} and the possibility to evolve the size and location of the local patch in order to capture the propagation of fracture using non-local damage phase-field equations in \cite{aldakheel_noii_wick_allix_wriggers_2021, aldakheel_noii_wick_wriggers_2021, noii_aldakheel_wriggers_2020, gerasimov_noii_allix_delorenzis_2018}.

Although a global-local approach relieves from highly time consuming modeling tasks (that are not scalable with computing resources), a chief drawback for the GLIC is the computational cost introduced with global-local iterations. The computing times of the GLIC are approximately equal to the total computing time of the direct simulation multiplied by the average number of global-local iterations during the simulation.
As summarized in \cite{Gosselet2018, blanchard_2019, Allix2022}, acceleration techniques can be employed as mitigation of increased computing costs and to improve robustness, similarly to fluid-structure interaction techniques as, for instance, the Anderson acceleration scheme \cite{anderson_1965}, more recently revisited in \cite{degroote_2009}.

With these works, the Authors explore and test the GLIC as implemented in Abaqus through co-simulation, which is a special technique that allows users to run simulations concurrently, exchanging quantity values through a finite element interface. In the case of the GLIC, co-simulation is run to couple the global simulation to the local one and to exchange displacements and reaction forces. The fields are exchanged through a user-defined surface in common between the two finite element models. A configuration file needs to be written to provide instructions to the Abaqus co-simulation engine on how to coordinate the two simulations.

Particular focus has been devoted on the study of the convergence and on the investigation of inexact convergence, as in classical multi-grid schemes as for instance \cite{fish_1992, fish_belsky_1995, fish_belsky_pandheeradi_1996, fish_suvorov_belsky_1997}, in order to reduce even further the computation expenses.
After this Introduction, Section~\ref{sec:formulation} is devoted to the formulation of the iterative scheme and the presentation of accelerators, Section~\ref{sec:cases} describes some challenging use cases and Section~\ref{sec:conclusions} proposes concluding remarks and future works.


\section{Formulation and study of the Global-Local Iterative Coupling}
\label{sec:formulation}
We consider a (quasi)-static mechanical problem set on a domain $\Omega$. We assume that a first \emph{global} (superscript $G$) finite element modeling is available, capable of correctly grasping the long-range trends in the structure. The global model need not be linear, in which case an incremental study is conducted, and for simplicity reason we only consider one load increment in this paper, multiple increments with time grids adapted to the models have been studied in~\cite{blanchard_2019}. The system to be solved can be written as:
\begin{equation}\label{eq:form}\text{Find }\bu\text{ so that }
	\foint\G(\bu\G)+\foext\G=0,
\end{equation}
where $\foext$ stands for the (generalized) external forces, $\foint$ is the vector of internal forces and $\bu$ is the vector of unknown displacements (Dirichlet conditions are assumed to have been eliminated).
Of course, in the case of a linear problem, the internal forces take the form of a linear application:
\begin{equation}
	\foint\G(\bu\G) = -\bK\G \bu\G
\end{equation}
where $\bK\G$ is the sparse symmetric definite positive stiffness matrix.

We suppose that this global model is insufficient in some \emph{local} areas: the geometry, the heterogeneity or the material law are too simplified, or the mesh is too coarse\ldots. These regions $(\Omega\Gs)_{s>0}$ are assumed to be non-overlapping sets of connected elements of the global mesh. For simplicity we consider that the meshes are matching at the interface, that is to say that there is a one-to-one correspondence between degrees of freedom. A possibility to implement this hypothesis is to have a copy of the global mesh of the interface on each patch and to use a tie constraint with the refined mesh.

\begin{figure}
    \centering
    \includegraphics{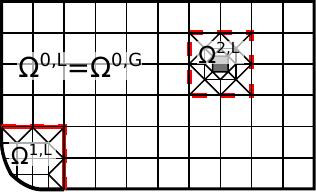}
    \caption{Reference problem with two local patches and a non-empty complement zone.}
    \label{fig:glref}
\end{figure}

On these regions, there exists a \emph{local} modeling (superscript $L$) with adapted geometry, mesh, elements, material law\ldots, nevertheless we suppose that the interface is preserved in the modelings: $\partial\Omega\Ls\cap\Omega=\partial\Omega\Gs\cap\Omega$. The zone of the global mesh which is not covered by patches, sometime referred to as complement zone, is written with index $0$: $\Omega\Gz=\Omega\G\setminus\cup_{s>0}\Omega\Gs$. Note that this zone may not exists if the global model is fully covered by patches. Our reference is the assembly of all the local models and of the complement zone: $\Omega\Gz\cup\left(\cup_{s>0}\Omega\Ls\right)$, see Figure~\ref{fig:glref}.

\begin{figure}
    \centering
    \includegraphics{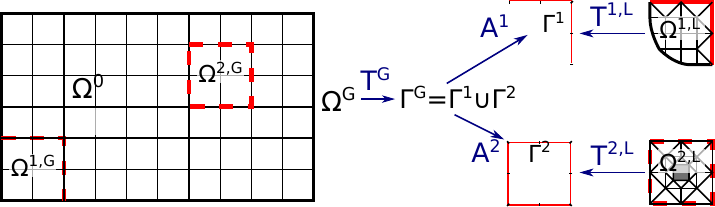}
    \caption{Global problem and topological operators (traces and assembly -- Use the transpose to follow the opposite direction).}
    \label{fig:glop}
\end{figure}

We note $\Gamma\s=\partial\Omega\s\cap\Omega$ the interface of the patches, and $\Gamma\G=\cup_{s>0}\Gamma\s$ the global interface. We introduce the trace operators $(\bT^G,\bT\Ls)$ that extract the interface degrees of freedom. We also introduce the injection operators $\bA\s$ which position $\Gamma\s$ with $\Gamma^G$, see Figure~\ref{fig:glop}.

We consider Dirichlet problems on the patches with the imposed displacement inherited from the global problem. Let Subscript $i$ denote \emph{internal} degrees of freedom of the patches and Subscript $\Gamma$ \emph{interface} degrees of freedom. 
\begin{equation}\label{eq:locdiri}
	\begin{aligned}
	&	\text{Dirichlet condition: }\bu_\Gamma\Ls = \bA\sT\bT\G\bu\G\\
	&	\text{Internal nonlinear problem: Find }\bu_i\Ls \text{ such that }  \fointi\Ls\begin{pmatrix}
			\bu\Ls_i \\ \bu_\Gamma\Ls
		\end{pmatrix} + \foexti\Ls = 0\\
	&\text{Interface post-processing: }\lam\Ls = - \fointb\Ls\begin{pmatrix}
		\bu\Ls_i \\ \bu_\Gamma\Ls
	\end{pmatrix} - \foextb\Ls
	\end{aligned}
\end{equation}
$\lam\Ls$ is the vector of nodal reaction associated with the imposed displacement $\bA\sT\bT\G\bu\G$, it can also be obtained as the Lagrange multiplier which imposes the given displacement.

At this point, we have carried out a classical sub-modeling approach. This approach is known to potentially brings lots of errors because the effect of the local models is not brought back on the global model, which in particular prevents interactions between patches.

The residual $\br_\Gamma$ can be measured as the lack of balance between the refined patches and the subdomain 0:
\begin{equation}\label{eq:resid}
\br_\Gamma = -\left(\bA\z\lam\Gz + \sum_{s=1}^N \bA\s \lam\Ls\right).
\end{equation}
The reaction on subdomain $0$ is defined as $\lam\Gz = - \fointb\Gz\begin{pmatrix}\bu\Gz_i \\ \bu_\Gamma\Gz \end{pmatrix} - \foextb\Gz$. A formula which avoids the management of Subdomain 0 will be given soon.

The global/local coupling simply consists in doing a Richardson iteration to try to reduce the residual.
To do so, we introduce the (initially null) corrective load $\bT\GT\bp\G_\Gamma$ applied to the global model.
The update of the load and its application in the corrected global problem can be written as:
\begin{equation}\label{eq:glocorr}
\begin{aligned}
        \bp\G_\Gamma \leftarrow \bp\G_\Gamma+ \omega \br_\Gamma,\\
	\foint\G(\bu\G)+\foext\G +\bT\GT\bp\G_\Gamma=0.
\end{aligned}
\end{equation}
where $\omega>0$ is a relaxation parameter. 

\begin{rem}
Other implementations of the global-local coupling, like the one in \cite{duval:hal-01065538}, prefer to use $\bu\G$ as the main unknown to be updated with relaxation at each iteration. When the global problem is linear, $\bu\G$ depends on $\bp_\Gamma$ in an affine manner and the two approaches are strictly equivalent. The approaches may differ slightly if the global problem is nonlinear, and we prefer our approach where the unknown has the same physical dimension (and vector size) as the ``natural'' residual which is the lack of balance between the models. Like in any other domain decomposition method, it highlights the crucial role played by the interface.
\end{rem}

The existence of a non-empty interval of values of the relaxation coefficient for which the algorithm converges is proved for an elastic global model coupled with monotonic patches (e.g. visco-plasticity with positive hardening). In practice the convergence is observed in more general cases, like with an elasto-visco-plastic global model or softening (but stable) patches. Moreover, all fixed-point accelerators can be used to improve the convergence.

Once the iteration has converged, we recover the \emph{reference} model constituted by the complement zone $\Omega\Gz$ and the patches $(\Omega\Ls)$ where the displacement is continuous at the interface and the nodal reaction are balanced, that is to say all interactions between the complement zone and the patches have been counted for. 

Algorithm~\ref{alg:staAit} recalls the classical global-local coupling method.
\begin{algorithm2e}[ht]\caption{Non-invasive stationary iterations with relaxation}\label{alg:staAit}\DontPrintSemicolon
\begingroup
\setlength{\belowdisplayskip}{0pt}
	Arbitrary initialization $\bp_{\Gamma,0}$\;
	\For{$j\in\left[0,\cdots,j_{\max}\right]$}{
     Global solve: Find $\bu\G_j$ such that:
     \[\foint\G(\bu\G_j)+\foext\G+\bT\GT\bp\G_{\Gamma,j}=0\]\;
     Global post-process of $\lam\Gz_j$ (if it exists):
     \[\lam\Gz_{j}=-\bT\Gz\left(\foint\Gz(\bu\Gz)+\foext\Gz\right)\]\;
     Local solves: $\forall s>0$, find $\bu\Ls_j,\lam\Ls_j$ such that: 
     \[
     \left\{\begin{aligned}
        &\bT\Ls\bu\Ls_{j} = \bA\sT\bT\G\bu\G_j\\
        &\foint\Ls(\bu\Ls_j)+\foext\Ls+\bT\LsT\lam\Ls_{j}=0
     \end{aligned}\right.
     \]\;
	Assembly of residual: $\br_{\Gamma,j}=  -\left(\bA\z\lam\Gz_j + \sum_{s=1}^N \bA\s \lam\Ls_j\right)$\;
 \lIf{$\|\br_{\Gamma,j}\|$ small enough }{Break}
 Update: $\bp_{\Gamma,j+1}=\bp_{\Gamma,j}+\omega\br_{\Gamma,j}$\;
	}
 \endgroup
\end{algorithm2e}

\begin{rem}
The post-processing of $\lam\Gz$ is sometimes difficult to implement in commercial software. A workaround consists in extracting the global version of the patches, solve the same Dirichlet problem as their refined counterparts and obtain the global nodal reactions on the patches $\lam\Gs$. Since the extra load $\bp_\Gamma$ introduces a stress discontinuity on the global model, we have:
\begin{equation}
    \sum_{s\geqslant 0}\bA\s \lam\Gs=\bp_\Gamma,
\end{equation}
and thus
\begin{equation}
\begin{aligned}
    \bA\z\lam\Gz &= \bp_\Gamma - \sum_{s> 0}\bA\s \lam\Gs\\
    \bp_\Gamma&\leftarrow (1-\omega)\bp_\Gamma + \omega\sum_{s> 0}\bA\s (\lam\Gs-\lam\Ls)
\end{aligned}
\end{equation}
\end{rem}

\subsection{Use of inexact solvers}
\label{sec:inexact}
In the case of nonlinear problems, the classical solution strategy is to use a Newton-Raphson or a quasi Newton iterative procedure. In that case, a stopping criterion must be used. Typically, the global solver can be written as: Find $\bu\G$ such that: 
\begin{equation}
     \foint\G(\bu\G)+\foext\G+\bT\GT\bp\G_{\Gamma}=\brho\G,\text{ with }\|\brho\G\|<\varepsilon\G_r.
\end{equation}

Another possible control is the use the length of the last correction to the unknown $\|\delta\bu\G\|<\varepsilon\G_\delta$. $(\varepsilon\G_r,\varepsilon\G_\delta)$ are the convergence thresholds (in term of residual and in term of variation). Those thresholds are often written in a relative way by introducing some normalization (e.g. the norm of the external load $\|\foext\G\|$). A classical question \cite{dembo1982inexact,eisenstat1996choosing}  is to try to tune the values of the thresholds in agreement with current outer loop residual, in our case $\|\br_\Gamma\|$, in order to avoid oversolving while preserving the good convergence rate. In this work, the authors make use of the Abaqus convergence controls~\cite{abaqus_2023}, which correpond to relative convergence criteria $\varepsilon\G_r=1$ (the default value being 0.005) and $\varepsilon\G_\delta=1$ (the default value being 0.01), where the normalizing factors are respectively the average flux norm and the norm of solution increment.

The fundamental equation to drive the thresholds, is the evaluation of the residual:
\begin{equation}
    \br_{\Gamma,j}=  -\left(\bA\z\lam\Gz_{j-1} + \sum_{s=1}^N \bA\s \lam\Ls_j\right) + \left(\bA\z\brho\G_{\Gamma,j}+\sum_{s>0}\bA\Ls\brho\Ls_{\Gamma,j}\right)
\end{equation}
In order for the global-local residual to be correctly evaluated, it is important to ensure that the residuals of the global and local solves remain negligible in comparison. If we use the notation $\Delta$ for the variation of quantities from one global-local iteration to another, we have:
\begin{equation}
    \bp_{\Gamma,j+1} =  \bp_{\Gamma,j} -\omega\left(\bA\z\Delta_j\lam\Gz + \sum_{s=1}^N \bA\s \Delta_j\lam\Ls\right) + \omega\!\left(\bA\z\brho\G_{\Gamma,j}+\sum_{s>0}\bA\Ls\brho\Ls_{\Gamma,j}\right)
\end{equation}
A good control appears to be $\|\brho\G_{\Gamma,j}\|\leqslant \alpha \|\Delta_j\lam\Gz\|$, with $\alpha\ll 1$, and the same relation for the patches. A simpler control is to use $\|\brho\G_{\Gamma,j}\|\leqslant \alpha \|\br_{\Gamma,j-1}\|$, and the same relation for the patches.


\subsection{Accelerators}
\label{sec:accelerators}
The Richardson iteration of Algorithm~\ref{alg:staAit} is the simplest example of fixed-point methods. Its convergence can be slow and it may even lack robustness, failing to converge on too stiff problems. This section investigates the use of 3 standard accelerators to overcome these issues.

For simpler presentation, we rephrase the coupling in terms of abstract operators:
\begin{equation}
\begin{aligned}
    \bu^G_\Gamma = \mathcal{S}^{G^{-1}}(\bp_\Gamma)& \qquad \text{obtained from eq.~}\eqref{eq:glocorr}\\
    \br_\Gamma = \mathcal{S}^L(\bu^G_\Gamma)& \qquad \text{obtained from eq.~}\eqref{eq:locdiri}\text{ and~}\eqref{eq:resid}
\end{aligned}
\end{equation}

The standard fixed-point system relies on the sequential Gauss-Seidel execution of solvers:
\begin{equation}
    \left(
        \mathcal{S}^L \circ \mathcal{S}^{G^{-1}}
    \right)
    \left(
        \bp_\Gamma
    \right)
    =
    \bp_\Gamma
\end{equation}
$S^G$ and $S^L$ refer to the computation of one global time increment and one local time increment respectively. 
%
%
%
The notations are further simplified by introducing the operator $H=\mathcal{S}^L \circ \mathcal{S}^{G^{-1}}$ and omitting the index $\Gamma$ in the unknown $\bp$:
\begin{equation}\label{FixedPointEquation}
    H: \mathbb{R}^n \to \mathbb{R}^n, H \left(\bp\right)=\bp
\end{equation}

In the reminder of this section, the lower $j$ denotes the iteration index to solve the fixed-point equation. $\tilde{\bp}_{j}$ denotes the $j-th$ Picard iterate  $\tilde{\bp}_{j}=H\left(\bp_{j}\right)$. For some problems, such as Fluid-Structure Interaction with moderate/large density ratios or structural multi-scale problems, a moderate to strong physical coupling bonds the subdomains. More elaborated algorithms such as the Aitken's relaxation or Quasi-Newton algorithms \cite{degroote_2010,uekermann_2016}, referred as accelerators, must be employed in order to get a stable solution and good convergence properties.
The cost associated with this additional acceleration step is negligible compared to the Picard iteration itself, especially for cases with surface coupling. Once the Picard iteration is computed, a correction is applied using accelerators:
\begin{equation}
    \bp_{j}
    \;\; \stackrel{\mathclap{\mbox{\tiny Picard}}}{\leadsto } \;\;
    \tilde{\bp}_{j}
    \;\; \stackrel{\mathclap{\mbox{\tiny Acc.}}}{\leadsto } \;\;
    \bp_{j+1}
\end{equation}
The classical residual can then be written as $\brx_j = \tilde{\bp}_{j} - \bp_{j}$.

\subsubsection{Aitken's relaxation}
The simplest accelerators are based on relaxation techniques:
\begin{equation}
    \bp_{j+1} = \bp_{j}\left(1-\omega_j\right)+\tilde{\bp}_{j}\omega_j
\end{equation}

A relaxation method with constant parameters might be sufficient for certain problems involving weak coupling. However, Aitken's relaxation is often used to improve the efficiency with a very limited additional computation cost. The Aitken's dynamic relaxation factor is provided by:
\begin{equation}
    \omega_j = \omega_{j-1}\frac{\brx_{j-1}^T\left(\brx_j-\brx_{j-1}\right)}{{\lVert \brx_j-\brx_{j-1}\rVert}_2^2}
\end{equation}


\subsubsection{Multi-secant methods}
In~\cite{fang_2009}, Fang and Saad present multi-secant methods to solve the nonlinear problem $\brx(\bp) = H(\bp)-\bp = 0$. The fixed-point problem is formally inverted in order to expose Picard iterations outputs:
\begin{equation}
    \brx \left(\bp\right)
    = 
    H\left(\bp\right)-\bp
    \approx
    \mathbf{0}
    \Rightarrow
    \tilde{\brx}\left(\tilde{\bp}\right)
    =
    \tilde{\bp} - H^{-1}\left(\tilde{\bp}\right) 
    \approx
    \mathbf{0}
\end{equation}

Introducing $\Delta \tilde{\bp}_{j} = \tilde{\bp}_{j+1} - \tilde{\bp}_{j}
$ and the Jacobian $J_{\tilde{\brx}} \left(\tilde{\bp}_{j}\right)$, the transition from the Picard step to the subsequent iterate in the Newton-Raphson update is articulated as follows:
\begin{equation}
    J_{\tilde{\brx}} \left(\tilde{\bp}_{j}\right)
    \Delta \tilde{\bp}_{j}
    =
    -\tilde{\brx}_{j}\left(\tilde{\bp}_{j}\right)
\end{equation}
\begin{equation}
    \bp_{j+1}
    =
    \tilde{\bp}_{j}
    +
    \Delta \tilde{\bp}_{j}
\end{equation}
Approximate of the interface Jacobian $J_{\tilde{\brx}}$ is computed using input-output information:
\begin{equation}
\begin{aligned}
    \bW_j &
    =
    \left[
        \Delta \tilde{\bp}_{0}, 
        \Delta \tilde{\bp}_{1}, 
        \ddots,
        \Delta \tilde{\bp}_{j-1}
    \right], &
    \text{with} \quad &
    \Delta \tilde{\bp}_{j} = \tilde{\bp}_{j+1} - \tilde{\bp}_{j}
    \\
    \bV_j &
    =
    \left[
        \Delta \brx_{0}, 
        \Delta \brx_{1}, 
        \ddots,
        \Delta \brx_{j-1}
    \right], &
    \text{with} \quad &
    \Delta \brx_{j} = \brx_{j+1} - \brx_{j}
\end{aligned}
\end{equation}

Using the two datasets $W_j$ and $Vj$, the multi-secant equation for the inverse interface Jacobian reads:
\begin{equation}
    J^{-1}_{\tilde{\brx}}\left( \tilde{\bp}_j \right) \bV_j \approx \bW_j
\end{equation}

\paragraph{Anderson acceleration \cite{anderson_1965}}
We seek for the minimum of the approximate inverse Jacobian norm while satisfying the multi-secant equation constraint. Let $L$ be the associated Lagrangian:
\begin{equation}
    L\left(J^{-1}_{\tilde{\brx}},\lambda\right):=\frac{1}{2} {\lVert J^{-1}_{\tilde{\brx}}\rVert}^2_F+\lambda \left(J^{-1}_{\tilde{\brx}} \bV_j-\bW_j\right) 
\end{equation}
$\lambda$ is the Lagrange multiplier associated with the multi-secant constraint. The solution $\left(\bar{\lambda}, \bar{J}^{-1}\right)$ is obtained by making stationary the Lagrangian. The inverse of the Jacobian can be written as:
\begin{equation}
    \bar{J}^{-1} = \bW_j \left( \bV_j^T \bV_j \right)^{-1} \bV_j^T
\end{equation}

The Anderson acceleration does not provide an approximation for the entire inverse Jacobian matrix; it focuses solely on the pertinent directions leading to the zero residual. We seek for the vector $\alpha$ minimizing the norm:
\begin{equation}
    \alpha_j = \arg\min_{\alpha} \lVert \bV_j \alpha +\brx_{j}\rVert_2
\end{equation}
The same direction is then applied to $\bW_j$ to get the prediction in terms of input variations:
\begin{equation}
    \bp_{j+1}=\tilde{\bp}_j+\bW_j \alpha_j
\end{equation}

Information from past increments can be used to improve the performance. However, the optimal number of reused time increments for the Anderson acceleration is highly dependent on the problem at stake. During the iterative process, dependent data within the data-sets $\bW_j$ and $\bV_j$ are removed using vector orthogonalization based on QR algorithms. 
Algorithm~\ref{alg:accAnderson} lists the details of the Anderson accelerator.
\begin{algorithm2e}[ht]
    \caption{Anderson acceleration}
    \label{alg:accAnderson}
    \DontPrintSemicolon
    \begingroup
    \setlength{\belowdisplayskip}{0pt}
    Extrapolation: $\bp_0=f(\bp^n, \bp^{n-1}, ...)$ \;
    Constant relaxation:
    $\tilde{\bp}_0=H(x_0)$, $\brx_0=\tilde{\bp}_0-\bp_0$ \;
    $\bp_1=\bp_0+\omega_0 \brx_0$ \;
	\For{$j\in\left[0,\cdots,j_{\max}\right]$}{
         Compute Picard iterate and residual: $\tilde{\bp}_j=H(\bp_j)$, $\brx_j=\tilde{\bp}_j-\bp_j$ \;
         \If{\normalfont{(!}$converged$\normalfont{)}}{
            Build $\bV_j$ and $\bW_j$ input-output datasets:
            $\bV_j= \left[ \Delta \brx_0, \dots, \Delta \brx_{j-1} \right]$, $\Delta \brx_i=\brx_{i+1}-\brx_i$ \;
            $\bW_j= \left[ \Delta \tilde{\bp}_0, \dots, \Delta \tilde{\bp}_{j-1} \right]$, $\Delta \tilde{\bp}_i=\tilde{\bp}_{i+1}-\tilde{\bp}_i$ \;
            Perform QR decomposition of $\bV_j$: $\bV_j=\mathbf{Q}\mathbf{U}$ \;
            Minimize $\lVert \bV_j \alpha +\brx_j\rVert_2$ using a least-square method: $\mathbf{U} \alpha = -\mathbf{Q}^T \brx_j$ \;
            Compute the correction and the accelerated solution:
            $\Delta \tilde{\bp}_j = \bW_j \alpha$ \;
            $\bp_{j+1} = \tilde{\bp}_j + \Delta \tilde{\bp}_j$ \;
         }
    }
    \endgroup
\end{algorithm2e}

\paragraph{Broyden's quasi-Newton method}
The Broyden's quasi-Newton method approximates the full inverse Jacobian matrix by minimizing the difference with the Jacobian at the previous increment under the multi-secant equation constraint:
\begin{equation}
    L\left(J^{-1}_{\tilde{\brx}},\lambda\right):=
    \frac{1}{2} {\lVert J^{-1}_{\tilde{\brx}} \left(\bp_j\right)-J^{-1, \left(N\right)}_{\tilde{\brx}} \left(\bp_j\right)\rVert}^2_F
    +\lambda \left(J^{-1}_{\tilde{\brx}} \bV_j-\bW_j\right) 
\end{equation}
Information from previous time increments are captured implicitly throughout this approach. The inverse of the Jacobian minimizing the Lagrangian function reads:
\begin{equation}
    \bar{J}^{-1} 
    = 
    J_{\tilde{\brx}}^{-1,\left(N\right)} 
    +
    \left(
        \bW_j - J_{\tilde{\brx}}^{-1,\left(N\right)} \bV_j
    \right)
    \left( \bV_j^T \bV_j \right)^{-1} \bV_j^T 
\end{equation}
Storing the complete inverse Jacobian matrix is impractical for real problems. Therefore, the inverse Jacobian storage is implemented through a truncated Singular Value Decomposition (SVD) to mitigate memory consumption. The SVD update is carried out using an incremental SVD \cite{scheufele_2017} to minimize computational costs. By using such an approach, the full rank SVD is never computed.

The solution correction $\Delta \tilde{\bp}_j$ is obtained using the Newton-Raphson update formula. The predicted (or the accelerated) value reads:
\begin{equation}
    \tilde{\bp}_{j+1} 
    = 
    \tilde{\bp}_{j} + \Delta \tilde{\bp}_j
    =
    \tilde{\bp}_{j} - J^{-1}_{\tilde{\brx}} \brx_j
\end{equation}

The pseudo-code of the Broyden's algorithm is detailed in Algorithm~\ref{alg:accBroyden}.
\begin{algorithm2e}[ht]
    \caption{Broyden's quasi-Newton method}
    \label{alg:accBroyden}
    \DontPrintSemicolon
    \begingroup
    \setlength{\belowdisplayskip}{0pt}
    Extrapolation: $\bp_0=f(\bp^n, \bp^{n-1}, ...)$ \;
    Constant relaxation: \;
    $\tilde{\bp}_0=H(x_0)$, $\brx_0=\tilde{\bp}_0-\bp_0$ \;
    $\bp_1=\bp_0+\omega_0 \brx_0$ \;
	\For{$j\in\left[0,\cdots,j_{\max}\right]$}{
         Compute Picard iterate and residual: $\tilde{\bp}_j=H(\bp_j)$, $\brx_j=\tilde{\bp}_j-\bp_j$ \;
         \If{\normalfont{(!}$converged$\normalfont{)}}{
            Build $\bV_j$ and $\bW_j$ input-output datasets: \;
            $\bV_j= \left[ \Delta \brx_0, \dots, \Delta \brx_{j-1} \right]$, $\Delta \brx_i=\brx_{i+1}-\brx_i$ \;
            $\bW_j= \left[ \Delta \tilde{\bp}_0, \dots, \Delta \tilde{\bp}_{j-1} \right]$, $\Delta \tilde{\bp}_i=\tilde{\bp}_{i+1}-\tilde{\bp}_i$ \;
            Perform QR decomposition of $\bV_j$: $\bV_j=\mathbf{Q}\mathbf{U}$ \;
            Compute $\mathbf{Z} =\left( \bV_j^T \bV_j \right)^{-1} \bV_j^T = \mathbf{U}^{-1}\mathbf{Q}^T$\;
            Update the approximate inverse Jacobian: \;
            $J_{\tilde{\brx}}^{-1} \brx_j \left(\tilde{\bp}_j\right)= 
            J_{\tilde{\brx}}^{-1,\left(n\right)} 
            +
            \left(
               \bW_j - J_{\tilde{\brx}}^{-1,\left(n\right)} \bV_j
            \right)
            \mathbf{Z}$ \;
            Compute the correction and the accelerated solution: \;
            $\Delta \tilde{\bp}_j = -J_{\tilde{\brx}}^{-1} \left(\tilde{\bp}_j\right)\brx_j$ \;
            $\bp_{j+1} = \tilde{\bp}_j + \Delta \tilde{\bp}_j$ \;
         }
    }
    \endgroup
\end{algorithm2e}

\begin{rem}
The GLIC is implemented in Abaqus through co-simulation. With this approach, the co-simulation coordinates two Abaqus concurrent simulations and their communications between each other. The quantities communicated in the GLIC case are displacements from global to local and reaction forces from local to global. Automatic mapping enables couplings between finite element interfaces with mismatching nodes, with a linear interpolation in space. If the simulations run with different increment sizes, a linear interpolation in time occurs. However, within this paper, only matching nodes and increment sizes have been considered. All the three accelerators described in Section~\ref{sec:accelerators} are implemented in the Abaqus co-simulation engine. From the practical point of view, beside running the global and local simulations concurrently, a so-called configuration file needs to be written and run in order to drive the Abaqus co-simulation engine. Within the configuration file, the user selects the accelerator technique and associated parameters, defines the field quantities that are exchanged and the pairs of outputs/inputs for, provides convergence and communication controls. Last but not least, it is important to note that users can substitute one of the two Abaqus simulations in the GLIC with their own third-party code, making use of special co-simulation APIs.
\end{rem}

\begin{rem}
It is important to state that, despite all the efforts to accelerate the convergence of the global-local iterations and optimize the number of global and local Newton-Raphson iterations, the GLIC will never be more efficient in terms of computing time than the monolithic finite elements model direct simulation, that will be taken as reference for accuracy measurements in what follows. More precisely, the GLIC will be slower by the factor close to the average number of global-local iterations during a given co-simulation run. However, for most of industrial cases, the computing time from simulations is often negligible in comparison to the time spent by engineers in building finite element models and, most importantly, is scalable with computing resources. The convenience of a global-local approach resides in the mitigation of the modeling burdens for creating multiple complex finite element models for multiple load cases or geometry configurations. Indeed, the GLIC proposes a unique global finite element model that is used for multiple simulations, whereas multiple specialized local models are applied where required, while keeping the global model unchanged. This is especially practical and a powerful approach in the aeronautic industry, where multiple levels of abstraction are defined for a given aircraft. The Authors refer to this concept as non-intrusive modeling.
\end{rem}

\section{Use cases}
\label{sec:cases}
A set of generic use-cases was defined to mimic situations met by Airbus where a direct unidirectional sub-modeling approach is not satisfactory:
\begin{itemize}
  \item{significant material non-linearity occurs locally in an area that is not modeled at global scale (plate with hole)}
  \item{significant material non-linearity occurs locally in an area that is modeled at global scale, but through a macro-element (elementary bolt)}
\end{itemize}

\subsection{Holed plate with localized elasto-plasticity}
\label{sec:case1}
\begin{figure}[htp]
    \centering
    \includegraphics[width=0.9\textwidth]{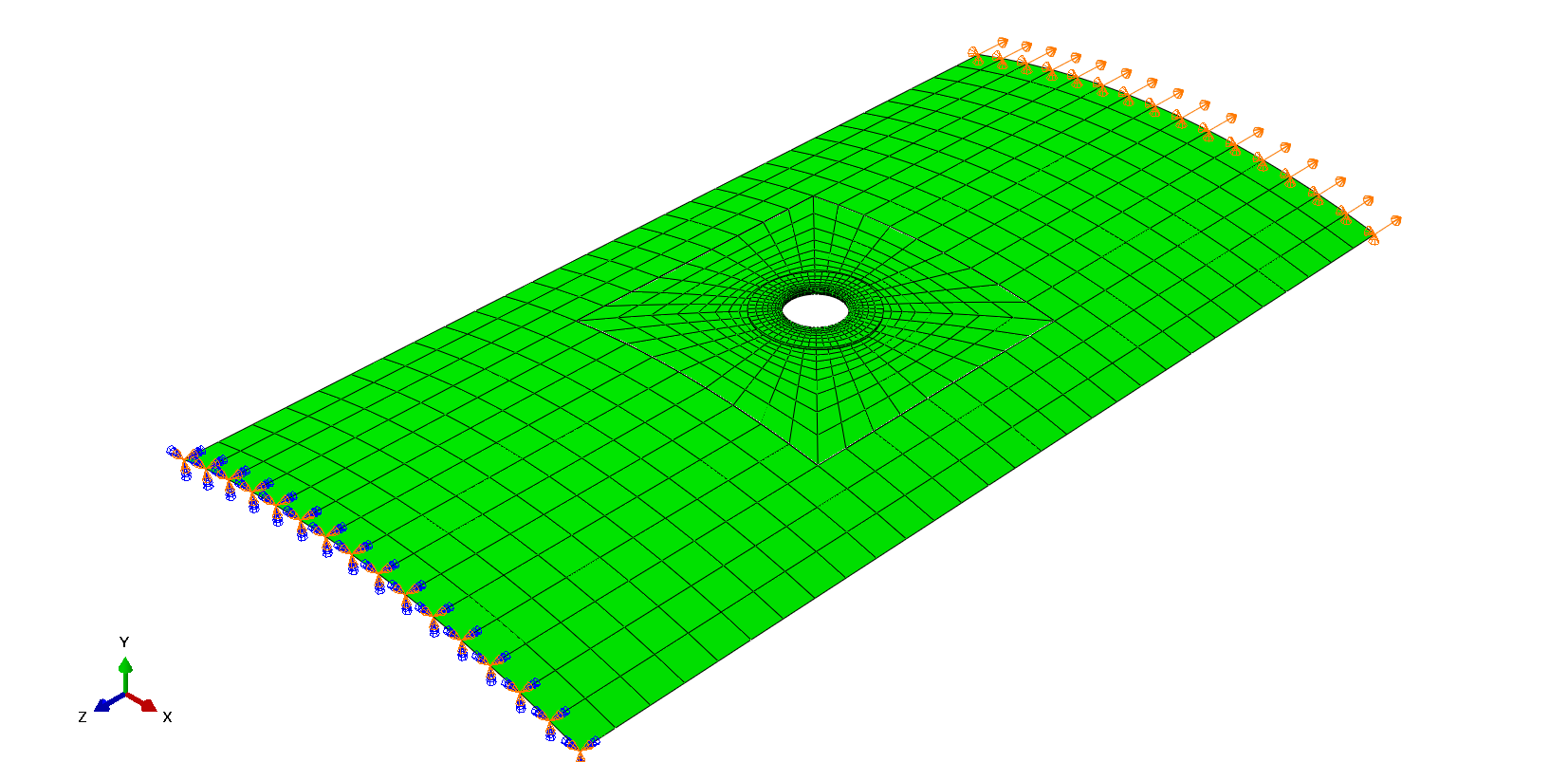}
    \caption{Holed plate finite element model of reference.}
    \label{fig:mesh1}
\end{figure}

The first case under study is a cylindrical holed plate section with localized elasto-plastic behavior. Figure~\ref{fig:mesh1} shows the monolithic finite element model of reference. The plate has width $10$ cm along the $X$ direction, length $20$ cm along the $Z$ direction, thickness $1$ mm and a cylindrical curvature along the $X$ direction of radius $0.2$ m. The hole is perfectly centered with respect to the plate and measures $6$ mm of radius. A mesh of solid elements is applied all around the hole to an outer radius of $12$ mm, with $8$-node hexahedral elements with reduced integration of variable size counting $8$ times through the thickness, $56$ times in the tangential direction through the full circle and $8$ times through the radial direction, for a total of $3584$ solid elements. The rest of the plate is meshed with $4$-node shell elements with Simpson reduced integration up to a mesh size of $6$ mm and connected to the solid elements with a shell-to-solid coupling (see~\cite{abaqus_2023} for details on this type of connection).

The plate is clamped at the bottom end at $Z=0.2$ m along the full width and pulled at the top end at $Z=0$ m with an imposed displacement of $1$ mm along the $Z$ direction. An elastic material is applied to the full model with a Young modulus of $210$ GPa and Poisson coefficient of $0.3$. The plastic material behavior is defined linearly interpolating the tabular values of Table~\ref{tab:plastic} and is applied to the solid elements around the hole. A nonlinear static analysis under the assumption of finite strains is therefore run with Abaqus, starting with a increment $50$ times smaller than the full step size and making use of the default convergence controls that adjust the increment size in function of convergence rates up to a maximum increment size $10$ times smaller than the full step size.

\begin{table}[ht]
\centering
\begin{tabular}{cc} 
 Stress (MPa) & Plastic strain \\
 \hline
 $400$ & $0.0$ \\ 
 $420$ & $0.02$ \\
 $500$ & $0.2$ \\
 $600$ & $0.5$ \\
 $625$ & $0.6$ \\
 $650$ & $0.8$ \\
\end{tabular}
\caption{Plastic material behavior on the solid elements around the hole.}
\label{tab:plastic}
\end{table}

\begin{figure}[htp]
    \centering
    \includegraphics[width=0.9\textwidth]{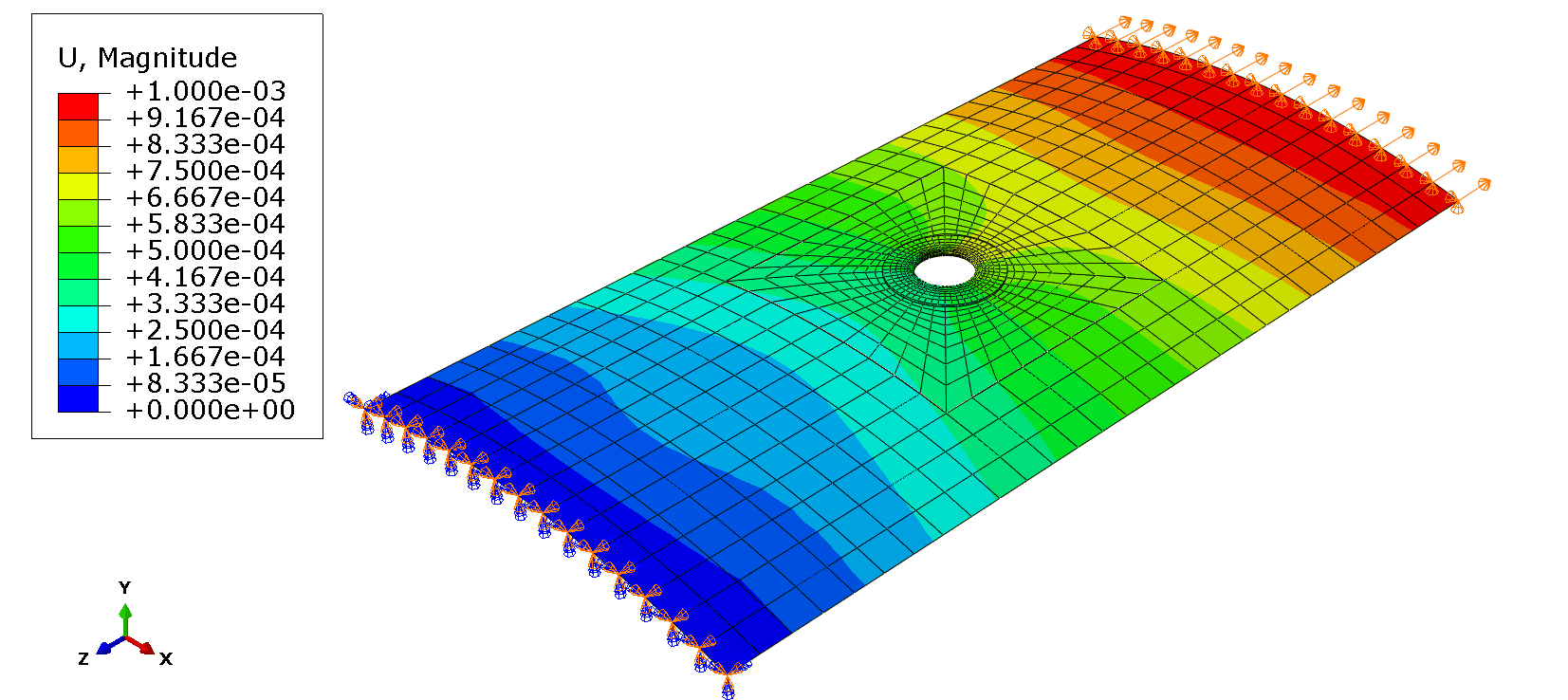}
    \caption{Results in terms of displacement magnitude (mm) for the holed plate finite element model of reference.}
    \label{fig:displ1}
\end{figure}

The contour plot in Figure~\ref{fig:displ1} shows results in terms of displacement magnitude projected to the deformed full mesh. All nodal displacements magnitudes vary between $0$ mm along the clamp and $1$ mm at the other end of the plate.

\begin{figure}[htp]
    \centering
    \includegraphics[width=0.9\textwidth]{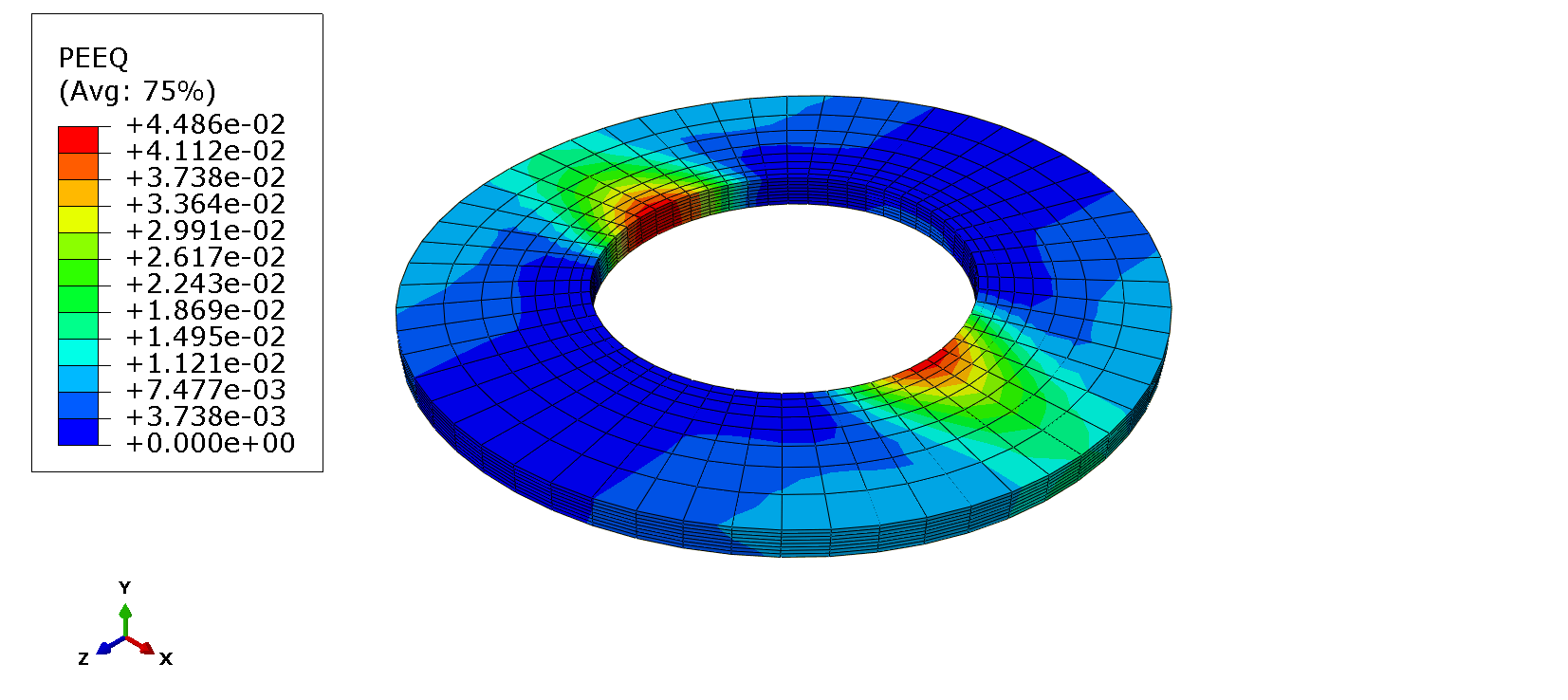}
    \caption{Results in terms of equivalent plastic strain for the holed plate finite element model taken as reference.}
    \label{fig:peeq1}
\end{figure}

The contour plot in Figure~\ref{fig:peeq1} shows results in terms of equivalent plastic strain projected to the solid elements mesh around the hole, which was the only part of the model subjected to elasto-plastic material behavior. Equivalent plastic strains reach a maximum plastic strain value of $0.04486$ on the side faces of the hole, where stresses concentrate.

The solution in Figures~\ref{fig:displ1} and~\ref{fig:peeq1} has been achieved with $26$ increments with a number of Newton-Raphson iterations varying between $2$ (at the beginning of the simulation) and $6$ (when yield stress is reached).

The GLIC approach is applied here with coupling a coarse global model with uniform mesh of quadrilateral shell elements of size $6$ mm without considering the presence of the hole to a local mesh that is identical to the mesh of the reference simulation on the squared patch of size $60$ mm around the hole. Automatic convergence controls of the global-local iterations have been implemented, namely, an increment-focused absolute criterion, an increment-focused relative criterion and a step-focused relative criterion, that have been set to $10^{-3}$ s, $10^{-3}$ s and $10^{-4}$ s, respectively, in what follows. This use case is challenging for the GLIC algorithm as the local model is different from the overlapped global region on multiple sides: element type, mesh size, material definition and evolution during the simulation. As the simple fixed-iteration coupling scheme converges too slowly for practical usage, only the accelerator techniques can be considered for this use case.

\begin{figure}[htp]
    \centering
    \begin{subfigure}{0.475\textwidth}
        \includegraphics[width=\textwidth]{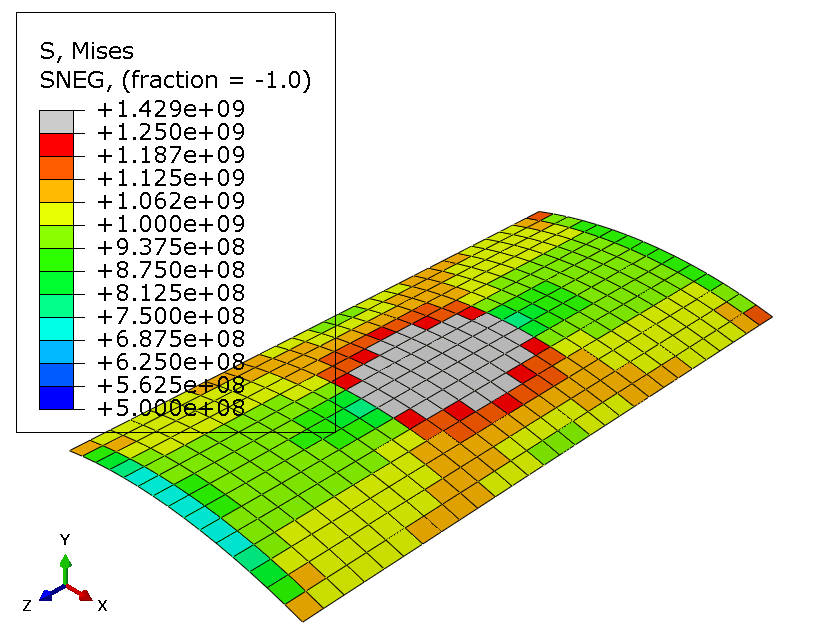}
        \caption{Global analysis.}
        \label{fig:glistressa}
    \end{subfigure}
    \begin{subfigure}{0.475\textwidth}
        \includegraphics[width=\textwidth]{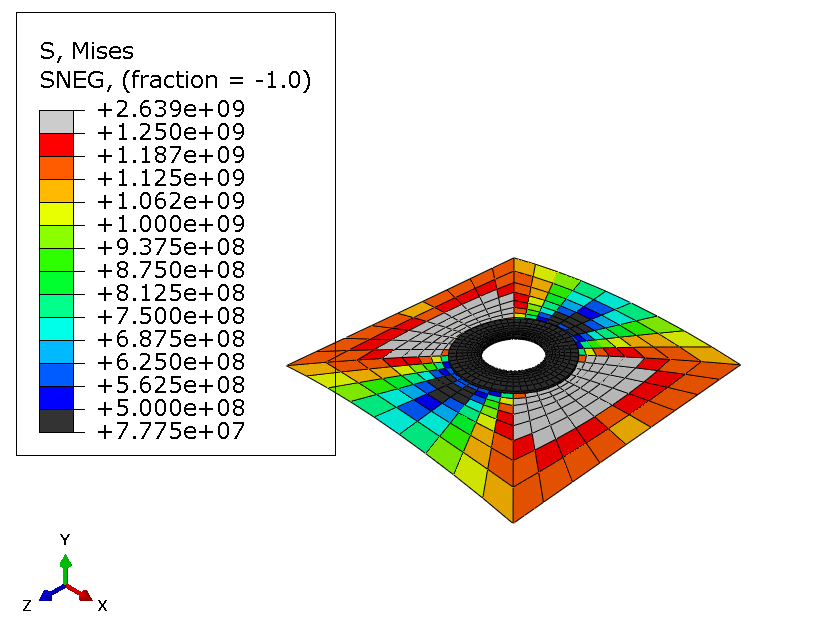}
        \caption{Local analysis.}
        \label{fig:glistressb}
    \end{subfigure}
    \caption{Results in terms of Von Mises stress (MPa) from the GLIC-Aitken simulation.}
    \label{fig:glistress}
\end{figure}

To appreciate the theoretical aspects of the GLIC, the contour plots from Figure~\ref{fig:glistress} show Von Mises stress fields in the global and local analyses for the co-simulation run with Aitken's relaxation (GLIC-Aitken). The contour plot color values vary in the range $[500,1250]$ MPa in order to qualitatively highlight the stress variation at the global-local interface. The overlapped region in the global model can be nicely identified with the larger jump in Von Mises stress values at the $9\times9$ elements squared patch over the center of the model: this is due to the correction forces in the global analysis that applied all along the global-local interface in the global model. On the other hand, the Von Mises stress values at the elements along the boundary of the local analysis are in accordance with the Von Mises stress values in the non-overlapped region of the global analysis.

\begin{figure}[p]
    \centering
    \begin{subfigure}{0.9\textwidth}
        \includegraphics[width=\textwidth]{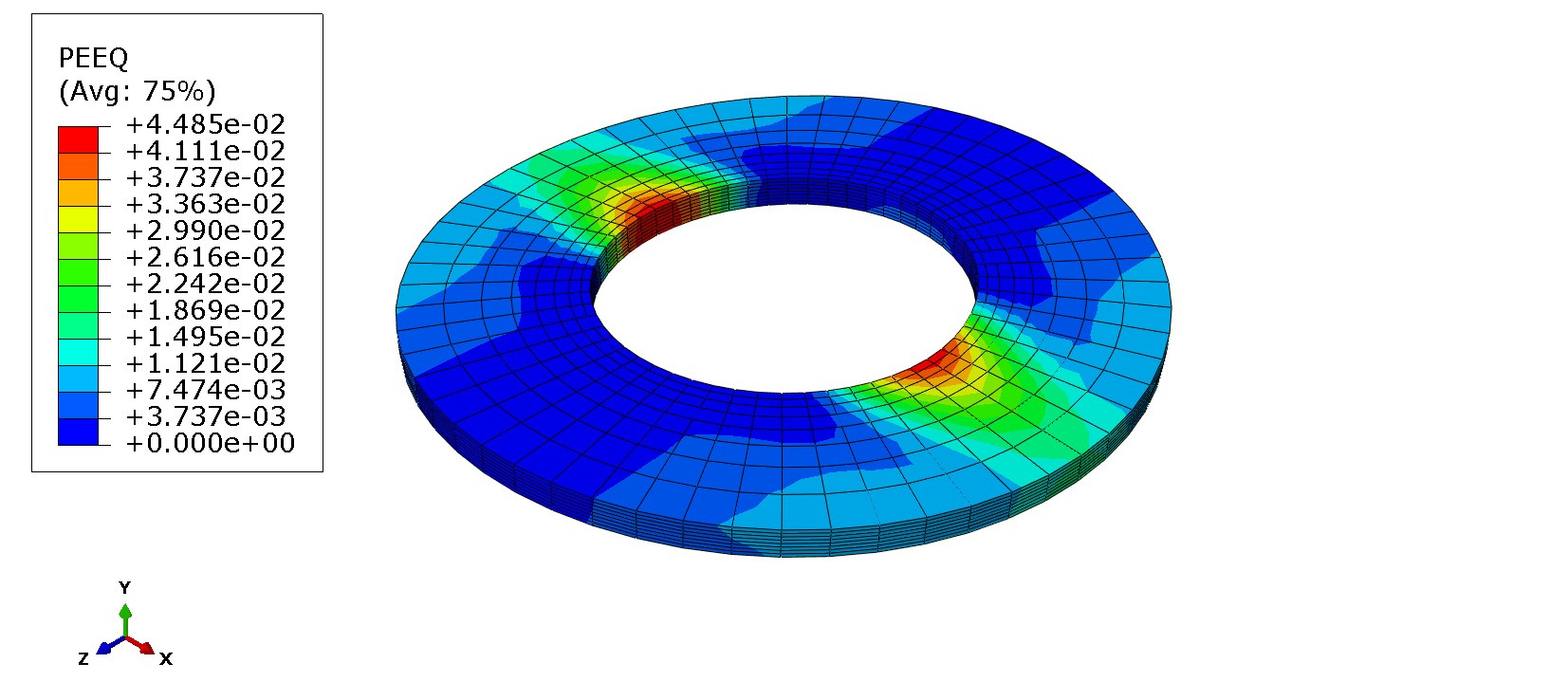}
        \caption{GLIC-Aitken simulation.}
        \label{fig:peeqacc1a}
    \end{subfigure}
    \begin{subfigure}{0.9\textwidth}
        \includegraphics[width=\textwidth]{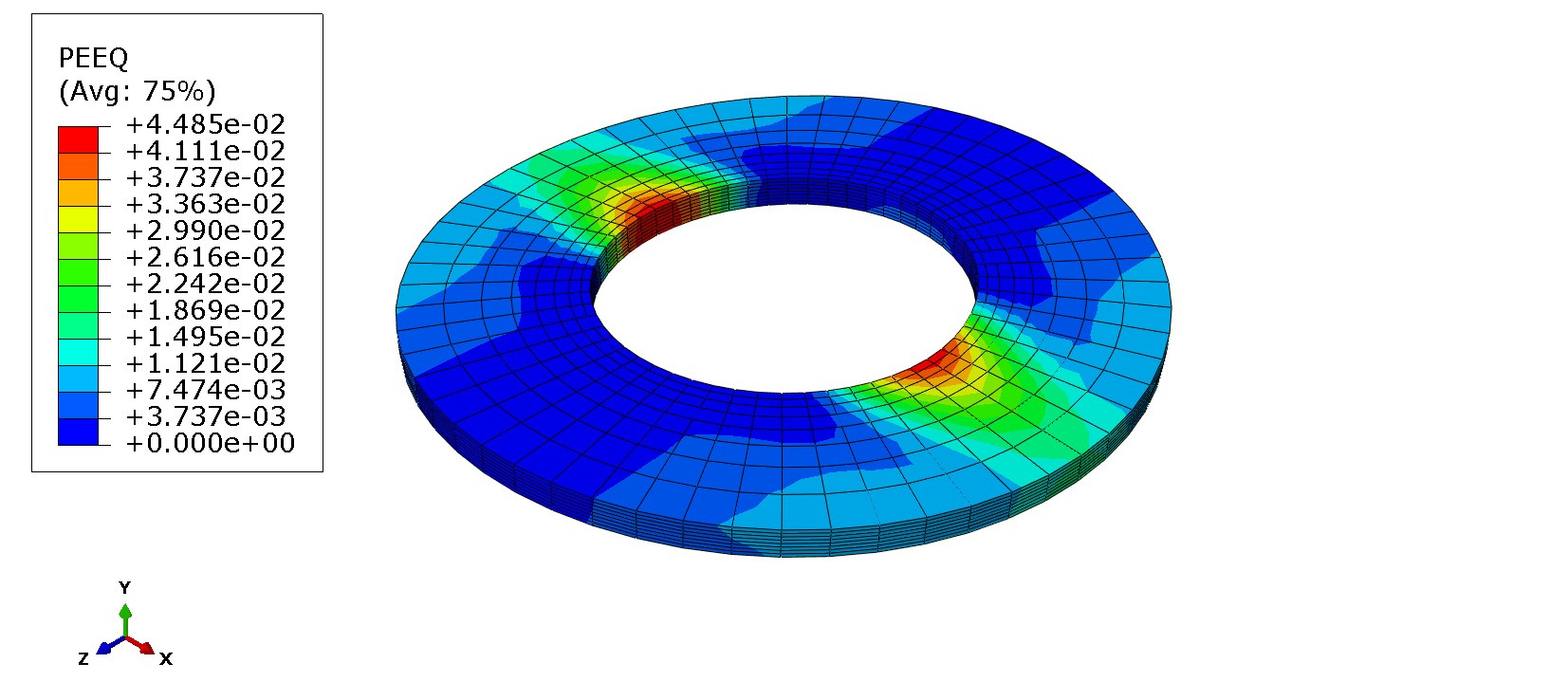}
        \caption{GLIC-Anderson simulation.}
        \label{fig:peeqacc1b}
    \end{subfigure}
    \begin{subfigure}{0.9\textwidth}
        \includegraphics[width=\textwidth]{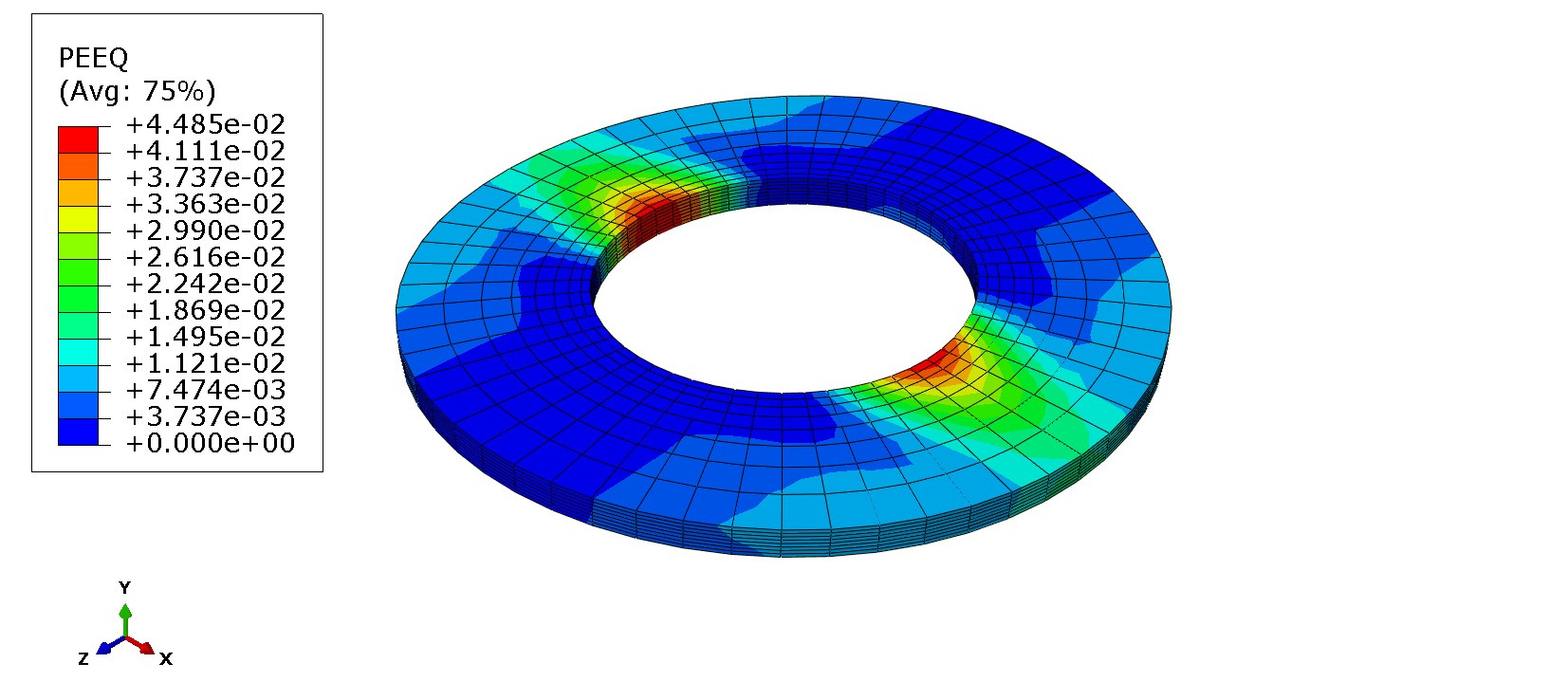}
        \caption{GLIC-Broyden simulation.}
        \label{fig:peeqacc1c}
    \end{subfigure}
    \caption{Results in terms of equivalent plastic strain for the GLIC simulations.}
    \label{fig:peeqacc1}
\end{figure}

In terms of accuracy, the contour plots in Figure~\ref{fig:peeqacc1} show results in terms of equivalent plastic strains obtained with Aitken's relaxation (GLIC-Aitken), Anderson acceleration technique (GLIC-Anderson) and Broyden quasi-Newton method (GLIC-Broyden) in Figures~\ref{fig:peeqacc1a},~\ref{fig:peeqacc1b} and~\ref{fig:peeqacc1c}, respectively. Such results are identical for engineering purposes and consistently show a maximum peak in equivalent plastic strain of $0.04485$, which is equivalent to an error of $0.02\%$ due to the same global-local convergence threshold for all three simulations.

\begin{figure}[p]
    \centering
    \includegraphics[width=0.74\textwidth]{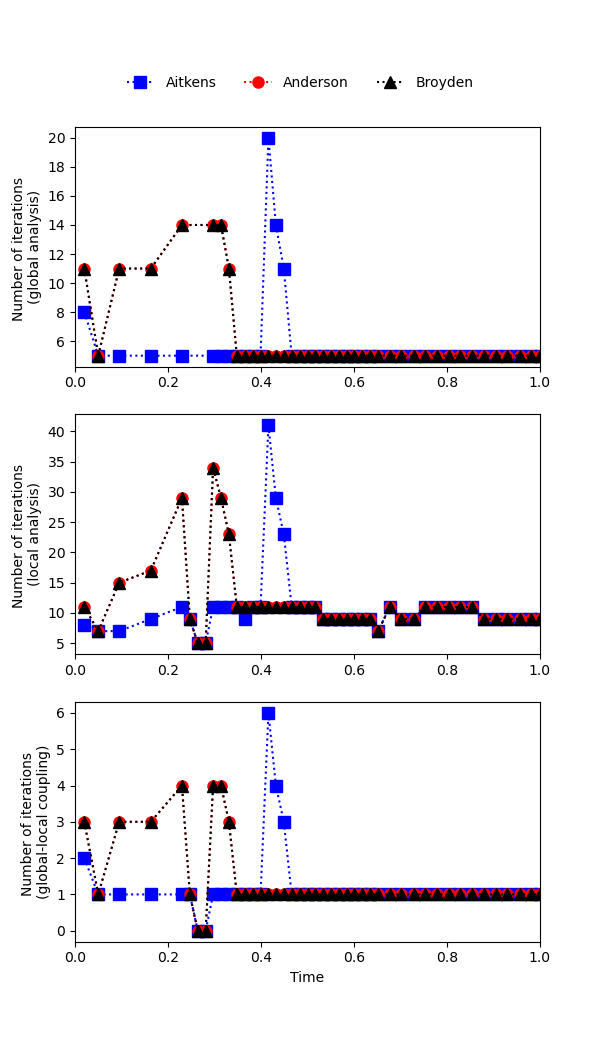}
    \caption{Comparison of the number of iterations in time between the GLIC-Aitken, GLIC-Anderson and GLIC-Broyden simulations.}
    \label{fig:iteracc1}
\end{figure}

In terms of performance, the numbers of iterations in time are reported in Figure~\ref{fig:iteracc1}. The plot on the top shows the total number of Newton-Raphson iterations in the global analysis (elastic finite strain model), the plot in the middle shows the total number of Newton-Raphson iterations in the local analysis (elasto-plastic finite strain model with coarse uniform mesh) and the plot on the bottom shows the number of iterations for the global-local coupling, which is a common multiplier of both global and local total Newton-Raphson iterations. The blue squares, the red circles and the black triangles are on curves out of the GLIC-Aitken simulation, GLIC-Anderson and GLIC-Broyden simulation, respectively. Whereas curves from the GLIC-Anderson and GLIC-Broyden simulations are perfectly overlapped with most of the iterations happening in the first half of the analysis, the GLIC-Aitken simulation curve behaves differently, with a sudden peak after time instant $0.4$ s, a few increments after a cutback in the increment size is triggered by the solver to optimize convergence.

The Authors want to note that the curves out of the global analysis quickly converge at each increment and global-local iteration, as the model is elastic. 

Overall, the GLIC-Aitken simulation is the fastest, independently from the sizes of global and local models. As highlighted in Table~\ref{tab:perfo1}, the total number of global-local iterations $N^{GL}$ is $1.5$ times smaller in GLIC-Aitken than in GLIC-Anderson and GLIC-Broyden. As a consequence, also in terms of total number of Newton-Raphson iterations in the global analysis $N^G_\text{iter}$ and in the local analysis $N^L_\text{iter}$, GLIC-Aitken is $1.57$ times faster than both GLIC-Anderson and GLIC-Broyden in the global analysis and $1.51$ times in the local analysis. Due to convergence control heuristics, the number of increments in the global analysis $N^G_\text{inc}$ and in the local analysis $N^L_\text{inc}$ are different. Namely, they are both smaller in the GLIC-Aitken than in both GLIC-Anderson and GLIC-Broyden, because of faster convergence on average. Last but not least, also the total number of global-local iterations is $1.52$ times smaller in GLIC-Aitken than in both GLIC-Anderson and GLIC-Broyden.

\begin{table}[ht]
\centering
\begin{tabular}{lccccc} 
  & $N^G_\text{inc}$ & $N^G_\text{iter}$ & $N^L_\text{inc}$ & $N^L_\text{iter}$ & $N^{GL}$ \\
 \hline
 GLIC-Aitken   & $32$ & $163$ & $35$ & $338$ & $65$ \\ 
 GLIC-Anderson & $41$ & $256$ & $44$ & $513$ & $99$ \\
 GLIC-Broyden  & $41$ & $256$ & $44$ & $513$ & $99$ \\
\end{tabular}
\caption{Comparison of total performance values between GLIC-Aitken, GLIC-Anderson and GLIC-Broyden.}
\label{tab:perfo1}
\end{table}

\begin{figure}[htp]
    \centering
    \includegraphics[width=0.9\textwidth]{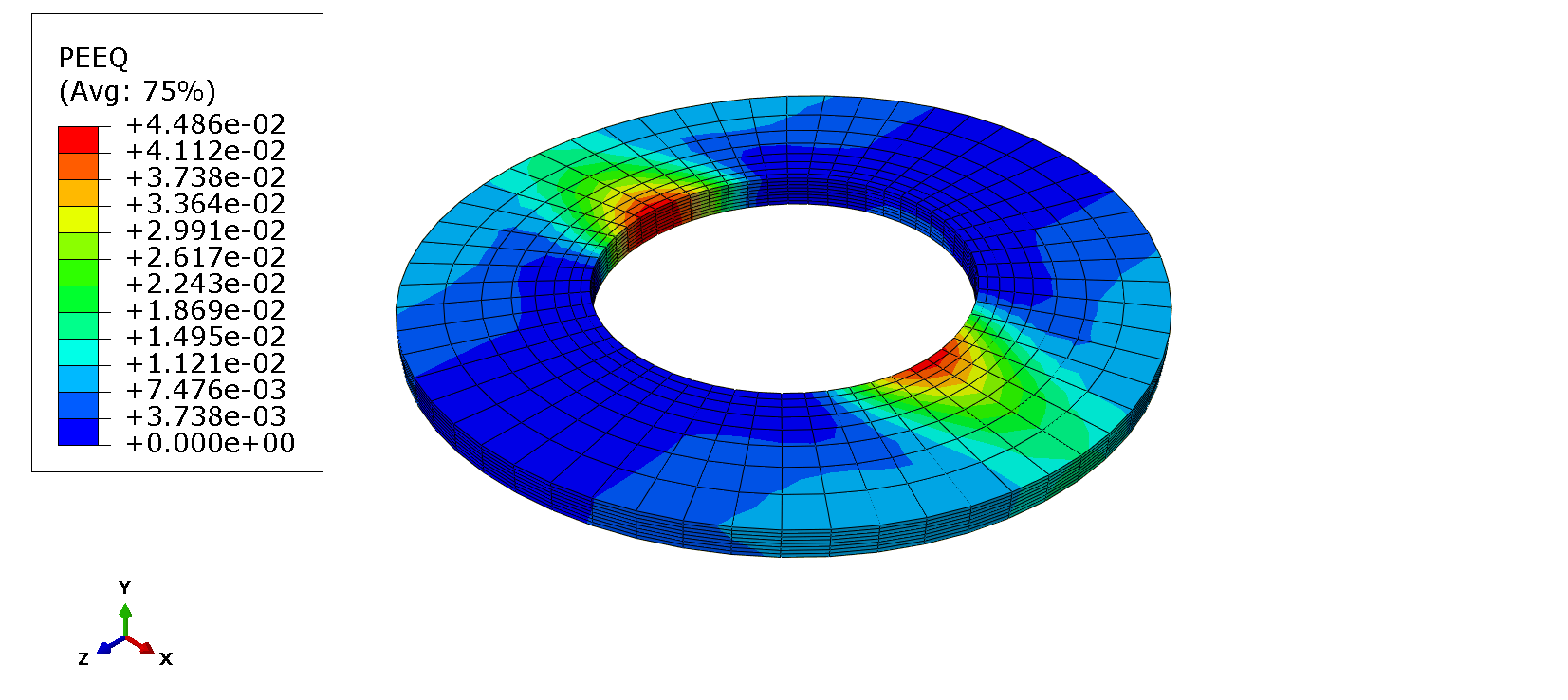}
    \caption{Results in terms of equivalent plastic strain in the local analysis from the GLIC-IS-Aiken simulation.}
    \label{fig:peeqaitkensctrls}
\end{figure}

The Newton-Raphson convergence controls employed in the simulation above have been set to the values that are proposed by default in the Abaqus solver. The readers are encouraged to consult the Abaqus manual~\cite{abaqus_2023} for comprehensive details on such convergence control values. To verify the inexact solver strategy introduced in Section~\ref{sec:inexact} (GLIC-IS), the following simulation have been run relaxing the convergence controls. Namely, the convergence criterion for the ratio of the largest residual to the corresponding average flux norm for convergence and the convergence criterion of the ratio of the largest solution correction to the largest corresponding incremental solution value have been increase from $0.005$ and $0.01$ to $1$ and $1$, respectively, for both global and local analyses. Contour plots in Figure~\ref{fig:peeqaitkensctrls} show the comparison in terms of equivalent plastic strain in the local model between the GLIC-Aitken and GLIC-IS-Aitken simulations.

\begin{figure}[htp]
    \centering
    \includegraphics[width=0.74\textwidth]{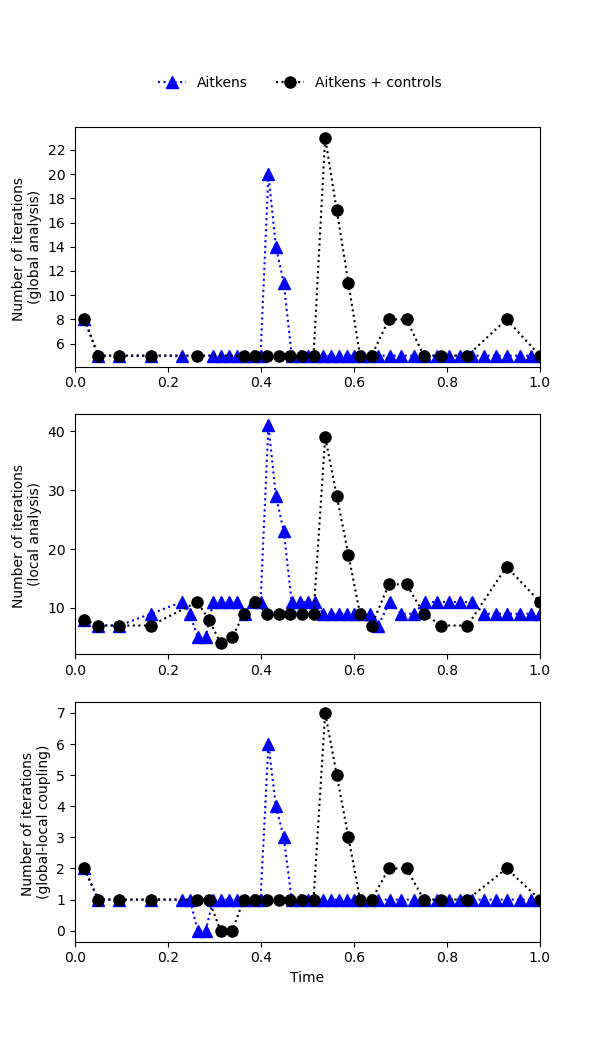}
    \caption{Comparison of the number of iterations in time between the GLIC-Aitken and GLIC-IS-Aitken simulations.}
    \label{fig:iterctrls1aitkens}
\end{figure}

Maximum values for the equivalent plastic strains are $0.04485$ for the GLIC-Aitken simulation and $0.04486$ for the GLIC-IS-Aitken simulation, with a negligible degradation in accuracy of approximately $0.02\%$. On the other hand, in terms of performance, the plots in Figure~\ref{fig:iterctrls1aitkens} show the comparison of the number of Newton-Raphson iterations during the global and local analyses between GLIC-Aitken and GLIC-IS-Aitken simulations. Although the curves have similar values in terms of iteration numbers on the vertical axes, the total number of iterations is lower in the relaxed case, due to larger increments, resulting in less equation solves. Indeed, the time increment is automatically increased or cut off based on convergence rate, which is faster with relaxed parameters.

As a more accurate performance analysis, the same numbers considered for the comparison between GLIC-Aitken, GLIC-Anderson and GLIC-Broyden are considered here in Table~\ref{tab:perfois1}, comparing GLIC-Aitken and GLIC-IS-Aitken. Whereas the total number of global-local iterations $N^{GL}$ is very similar in both simulations, as well as the total number of Newton-Raphson iterations in the global analysis $N^G_\text{iter}$, which is counterintuitively $3\%$ smaller in GLIC-Aitken than in GLIC-IS-Aitken. The total number of Newton-Raphson iterations in the local analysis $N^L_\text{iter}$ is $11\%$ smaller in GLIC-IS-Aitken than in GLIC-Aitken, due to small numbers of increments in both global and local analyses, $N^G_\text{iter}$ and $N^L_\text{iter}$, respectively.

\begin{table}[ht]
\centering
\begin{tabular}{lccccc} 
  & $N^G_\text{inc}$ & $N^G_\text{iter}$ & $N^L_\text{inc}$ & $N^L_\text{iter}$ & $N^{GL}$ \\
 \hline
 GLIC-Aitken    & $32$ & $163$ & $35$ & $338$ & $65$ \\ 
 GLIC-IS-Aitken & $24$ & $168$ & $27$ & $304$ & $64$ \\
\end{tabular}
\caption{Comparison of total performance numbers between GLIC-Aitken and GLIC-IS-Aitken.}
\label{tab:perfois1}
\end{table}

This time, the convenience in terms of performance depends on the size of the global and local models. As the number of variables in the local model is $16269$, much larger than the number of variables in the global model, that is $2880$, then in this case GLIC-IS-Aitken is more convenient than GLIC-Aitken.

Similar behaviors are encountered relaxing convergence criterion values with the Anderson acceleration and Broyden's method.

\subsection{Structure with bolted joints}
\label{sec:case2}
\begin{figure}[htp]
    \centering
    \includegraphics[width=0.9\textwidth]{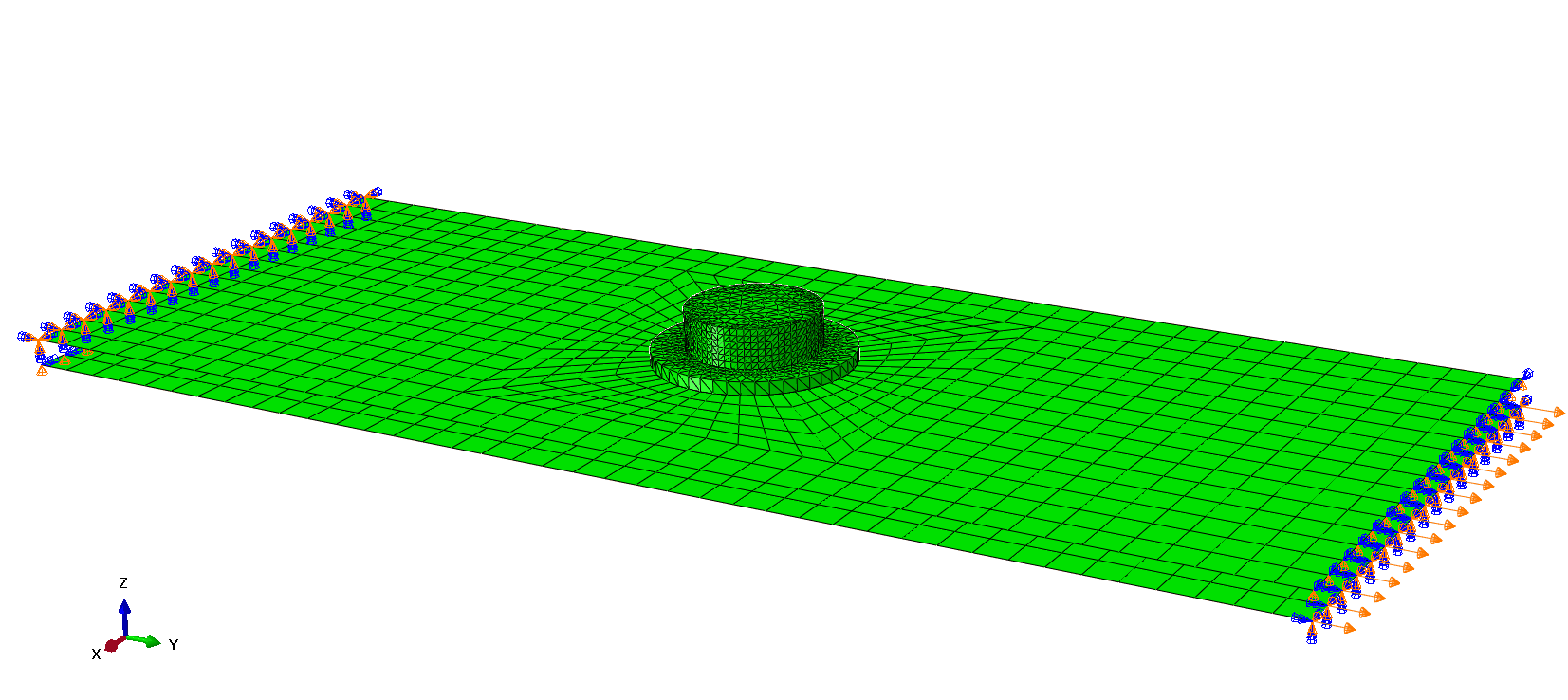}
    \caption{Bolted joint finite element model of reference.}
    \label{fig:pwbdirectmodel}
\end{figure}

The second case under study is a bolted joint connecting two planar plates together. Even if derived into a simplistic use-case (Figure~\ref{fig:pwbdirectmodel}), the situation explored here is representative of actual industrial concerns where elementary bolts are represented by connector elements at global scale and may meet critical loading conditions, which requires further refined investigations with local models of higher physical representability. The full model is an assembly composed of $3$ parts in contact between each other: the two plates and the bolt. Figure~\ref{fig:pwbdirectmodel} shows the reference finite element model applied to the geometry. The two plates are identical, of width $40$ mm along the $X$ direction, length $100$ mm along the $Y$ direction, thickness $2$ mm along the $Z$ direction and are stacked in the $Z$ direction. As for the holed plate, the bolt is perfectly centered across the two plates, its head has a diameter $10$ mm and height $3$ mm, its body has a diameter $6.35$ mm and a height $7$ mm. A mesh of solid elements is applied around the bolt up to a diameter $15$ mm, with $4$-noded tetrahedric elements of variable size, with smaller ones over the surfaces in contact. As for the holed plate, the rest of the assembly is meshed with $4$-noded shell elements with Simpson reduced integration up to a mesh size $2.5$ mm.

As one of the plates is clamped on the bottom edge and simply supported in the $Z$ direction on the top edge, the other plate is simply supported in the $Z$ direction on the bottom edge and pulled up in the $Y$ direction on the top edge for a displacement of $2$. Elastic materials are applied to the plates and the bolt with Young moduli $70$ GPa and $110$ GPa, respectively, and a Poisson coefficient $0.3$ in both cases. Small sliding contact with friction $0.2$ and slip tolerance $0.005$ mm is defined between the plates and between the bolt and the plates. A nonlinear static analysis under the assumption finite strains is run with Abaqus in two steps. In the first step, a pre-load of $1$ kN is applied to the bolt as compression in the $Z$ direction through bushing connectors of elasticity $100$ kN$/$mm for axial forces against translations and $10^5$ for moments against rotations. In the second step, one of the plate is pulled and loads the bolt with shear stress through contact.

\begin{figure}[htp]
    \centering
    \includegraphics[width=0.9\textwidth]{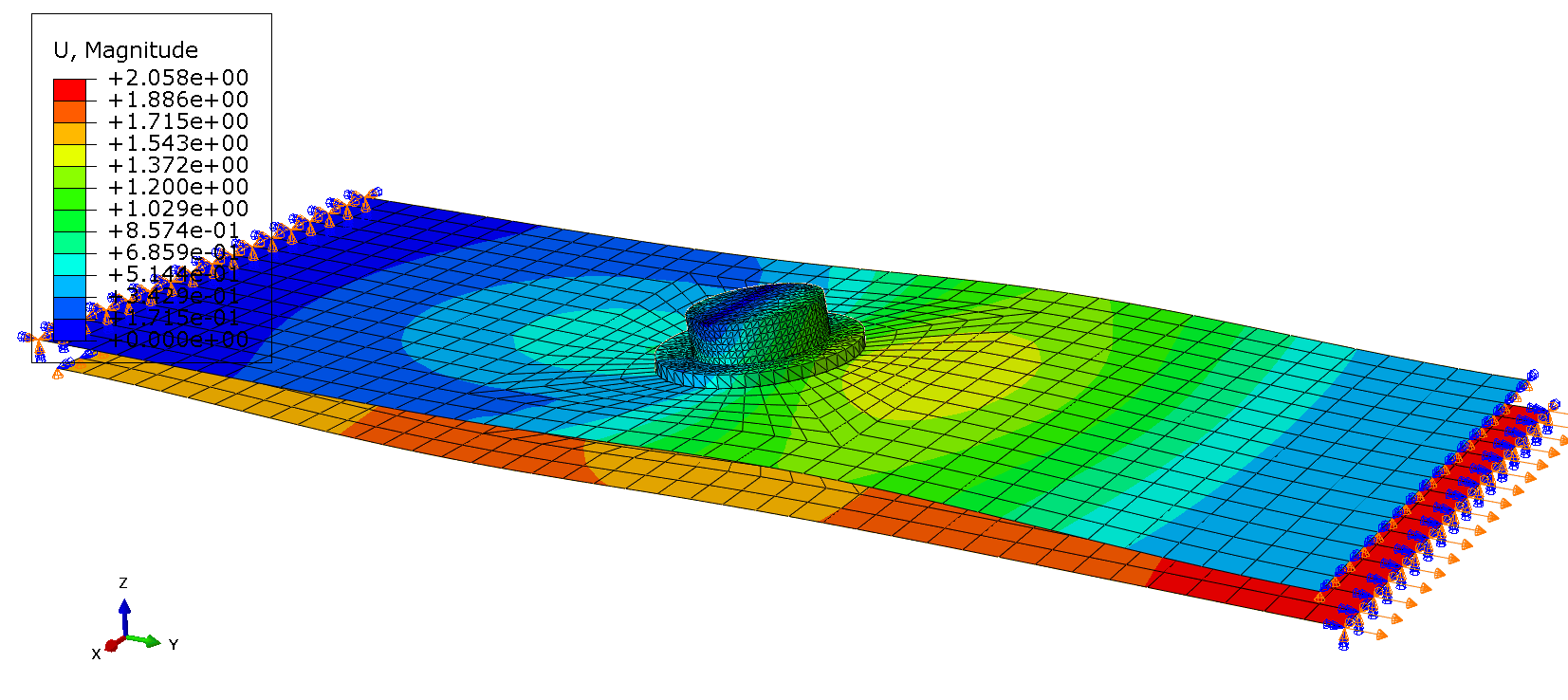}
    \caption{Results in terms of displacement magnitude (mm) for the bolted joint assembly finite element model of reference.}
    \label{fig:pwbdispl1}
\end{figure}

The contour plot in Figure~\ref{fig:pwbdispl1} shows results in terms of displacement magnitude projected to the deformed full mesh. All nodal displacements magnitudes vary between $0$ mm along the clamp and $2.058$ mm.

\begin{figure}[htp]
    \centering
    \begin{subfigure}{0.475\textwidth}
        \includegraphics[width=\textwidth]{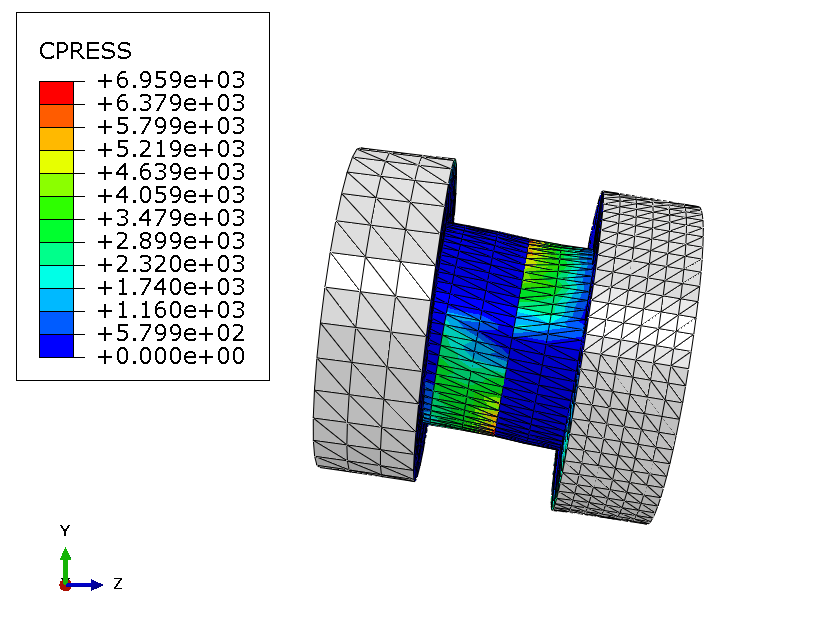}
        \caption{Contact pressure.}
    \end{subfigure}
    \begin{subfigure}{0.475\textwidth}
        \includegraphics[width=\textwidth]{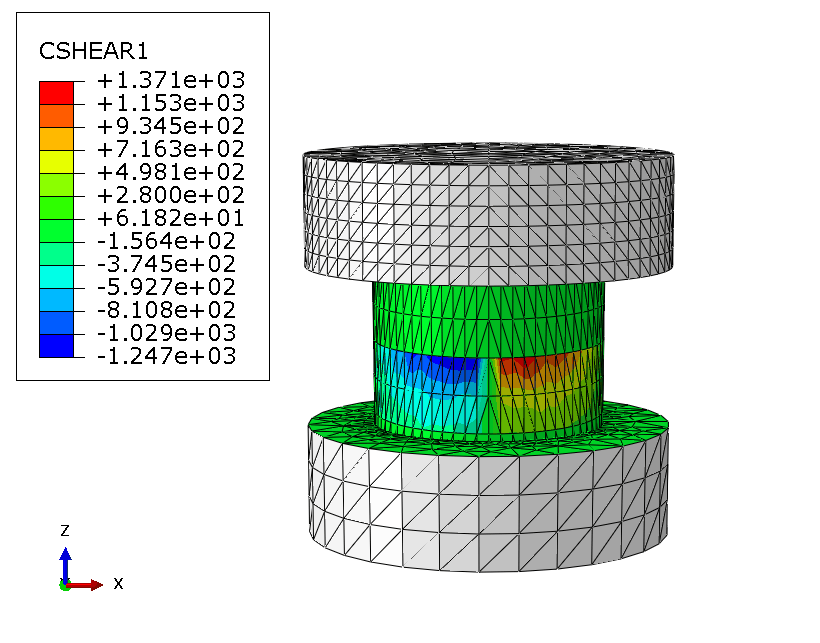}
        \caption{Contact shear in $Z$ direction.}
    \end{subfigure}
    \caption{Results in terms of contact interaction (MPa) on the bolt finite element model taken as reference.}
    \label{fig:pwbcontact}
\end{figure}

The contour plot in Figure~\ref{fig:pwbcontact} shows results in terms of contact pressure and contact shear, respectively, projected to the solid elements mesh of the bolt. Contact pressure goes up to an absolute value of $6.96\cdot10^3$ MPa in the middle of the bolt body in an almost symmetric manner. Contact shear participates with values ranging between $-1.25\cdot10^3$ MPa and $1.37\cdot10^3$ MPa.

The solution in this case has been achieved with 15 increments per step. The first step ran with a number of Newton-Raphson iterations varying between $1$ and $7$, on average $2.1\overline{3}$, the peak of $7$ caused by contact severe discontinuities. The second step ran with a number of Newton-Raphson iterations varying between $5$ and $9$, on average $6.2\overline{6}$, in part due to contact, in part due to geometrical nonlinearities.

In accordance to the GLIC strategy and similarly to the first use case described in Section~\ref{sec:case1}, the global model is meshed with a uniform coarse mesh of $4$-noded shell elements of mesh size $2.5$ mm, whereas the local model is defined on the central $30\times30$ mm$^2$ squared patch that covers the bolt and the immediate area of the two plates, meshed partly with $4$-noded shell elements and partly with $3$-noded tetrahedric elements, connected with shell-to-solid coupling. The overlap of the local model to the global model exactly represents what was defined in the reference model of Figure~\ref{fig:pwbdirectmodel}, as the goal of this verification exercise is to obtain identical results to reference. Automatic convergence control tolerances of $10^{-4}$, $10^{-3}$ and $10^{-3}$ have been used to stop global-local iterations as increment-focused absolute criterion, increment-focused relative criterion and step-focused relative criterion, respectively.

As the bolt is fully modeled in 3D in the local model, its behavior is captured in the global model with a bushing connector with elastic behavior of $275$ kN/mm in $X$ and $Y$ directions and $348\cdot 10^3$ kN/mm in the $Z$ direction.

This case is particularly challenging for the fixed-point GLIC algorithm as the local model is much more refined and globally stiffer than the global overlapped region and the spectral radius of the amplification matrix might likely be larger than $1$, which would lead to divergence. Indeed, also in this case, the fixed-iteration GLIC algorithm without any acceleration failed to converge. Therefore, only the results obtained with accelerator techniques are detailed in what follows.

\begin{figure}[htp]
    \centering
    \begin{subfigure}{0.475\textwidth}
        \includegraphics[width=\textwidth]{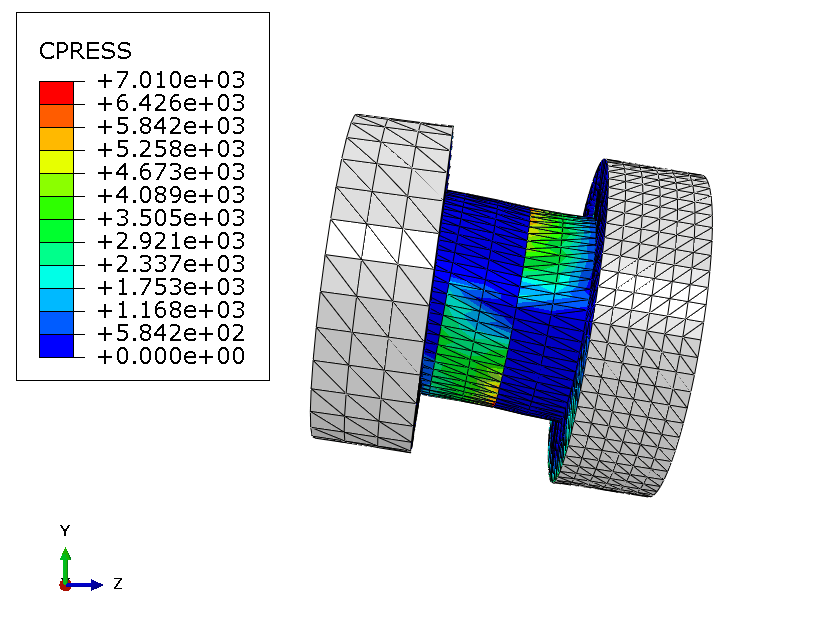}
        \caption{Contact pressure with GLIC-Aitken.}
    \end{subfigure}
    \begin{subfigure}{0.475\textwidth}
        \includegraphics[width=\textwidth]{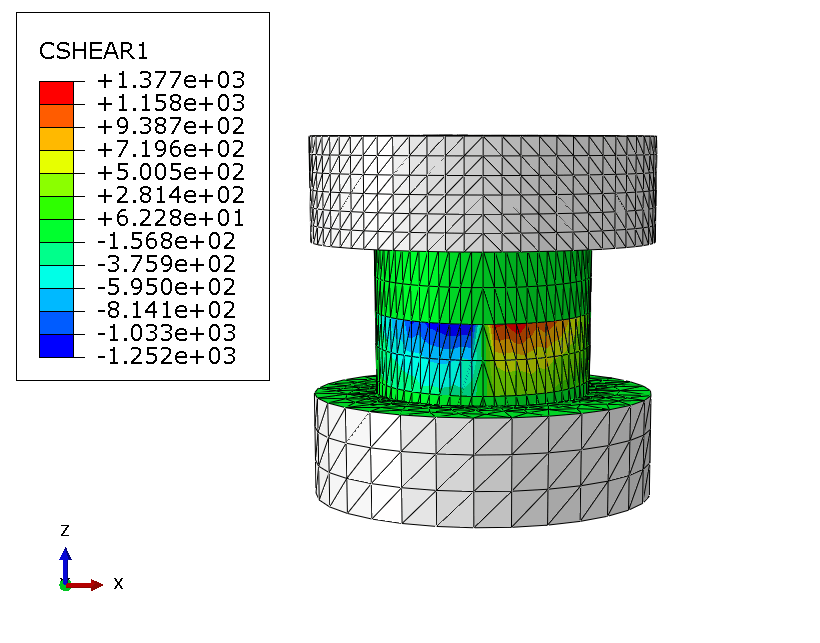}
        \caption{Contact shear ($Z$ dir.) with GLIC-Aitken.}
    \end{subfigure}
    
    \begin{subfigure}{0.475\textwidth}
        \includegraphics[width=\textwidth]{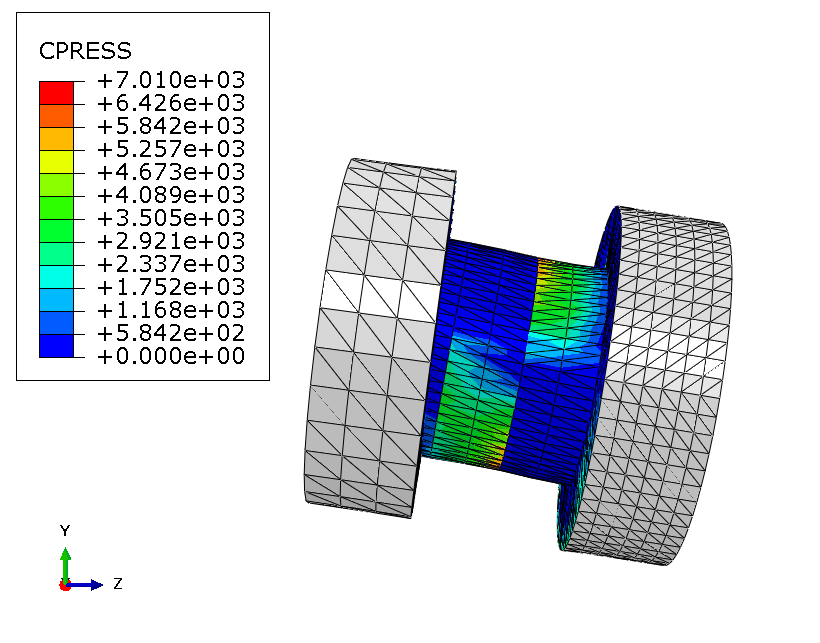}
        \caption{Contact pressure with GLIC-Anderson.}
    \end{subfigure}
    \begin{subfigure}{0.474\textwidth}
        \includegraphics[width=.98\textwidth]{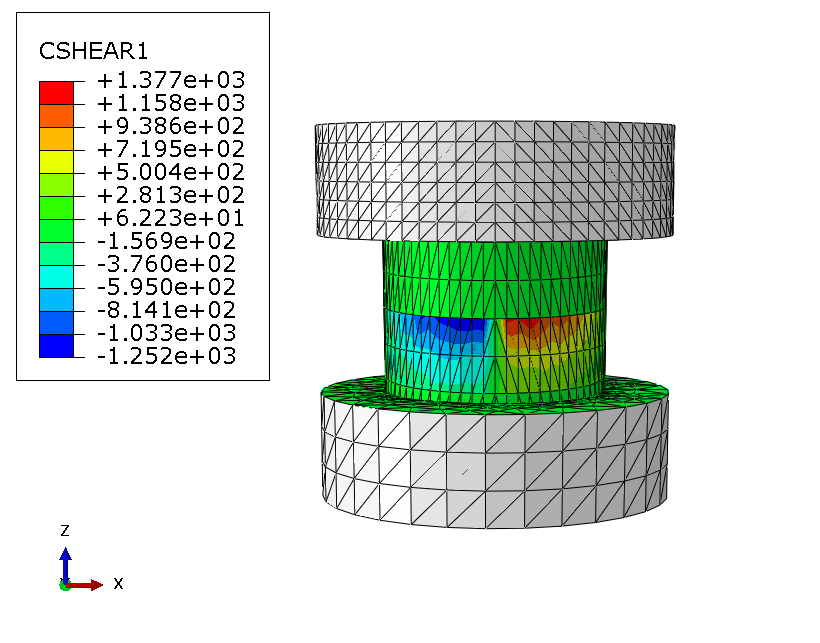}
        \caption{Contact shear ($Z$ dir.) with GLIC-Anderson.}
    \end{subfigure}
    
    \begin{subfigure}{0.475\textwidth}
        \includegraphics[width=\textwidth]{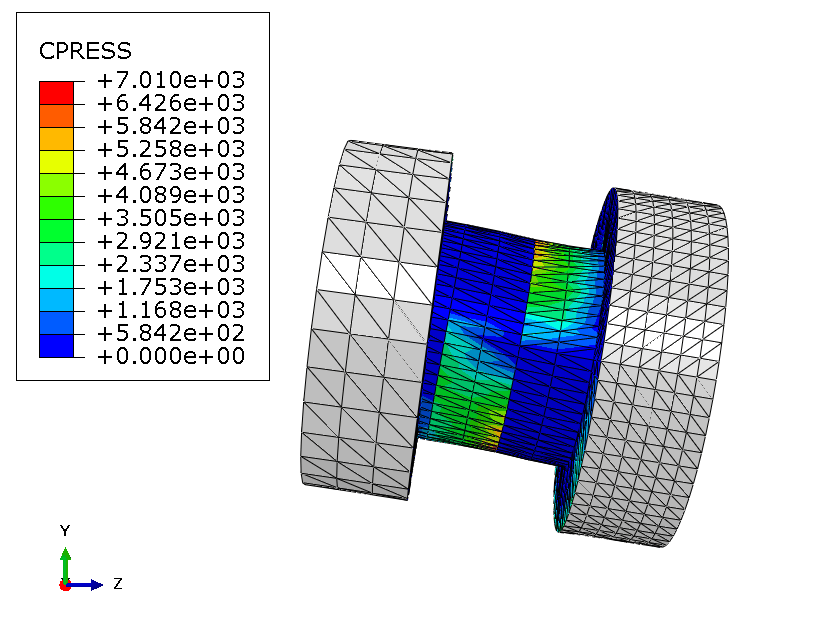}
        \caption{Contact pressure with GLIC-Broyden.}
    \end{subfigure}
    \begin{subfigure}{0.474\textwidth}
        \includegraphics[width=.98\textwidth]{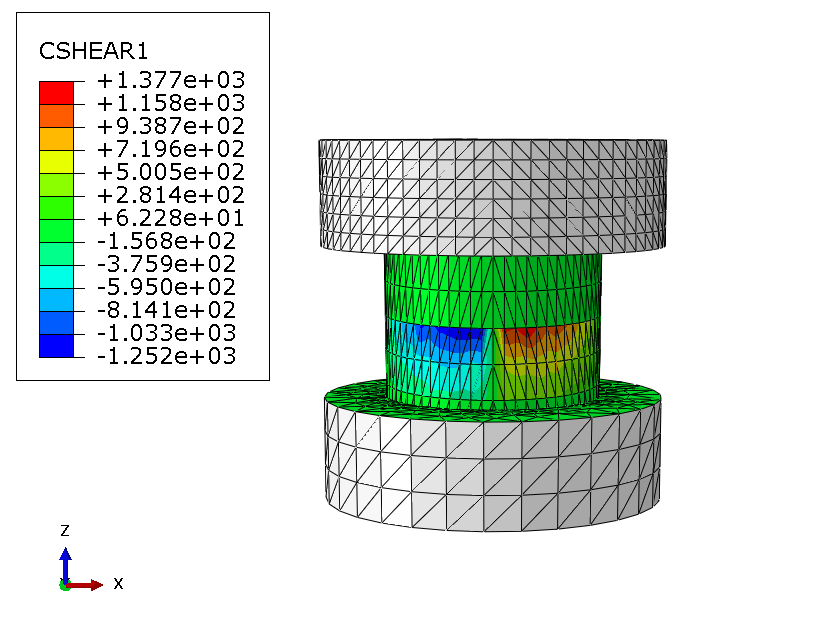}
        \caption{Contact shear ($Z$ dir.) with GLIC-Broyden.}
    \end{subfigure}
    \caption{Results in terms of contact interaction (MPa) on the bolt for the GLIC simulations.}
    \label{fig:pwbcontactacc}
\end{figure}

In terms of accuracy, Figure~\ref{fig:pwbcontactacc} shows contour plot results in terms of contact pressure and contact shear in the $Z$ direction for co-simulation runs with GLIC-Aitken, Anderson accelerator technique and GLIC-Broyden. Results are qualitatively identical between each other and are close enough to the reference solution, as detailed on maximum and minimum values with Table~\ref{tab:usecase2errors}. This possibly indicates that the error is chiefly due to the GLIC convergence criteria rather than the chosen accelerator methodology.

\begin{table}[htp]
\centering
\begin{tabular}{lccc} 
 & Max CPRESS & Max CSHEAR1 & Min CSHEAR1 \\
 \hline
 Reference     & $6.959$ & $1.371$ & $-1.247$\\ 
 GLIC-Aitken   & $7.01$  & $1.377$ & $-1.252$\\
 GLIC-Anderson & $7.01$  & $1.377$ & $-1.252$\\
 GLIC-Broyden  & $7.01$  & $1.377$ & $-1.252$\\
\end{tabular}
\caption{Contact interaction results comparison (values expressed in GPa).}
\label{tab:usecase2errors}
\end{table}

\begin{figure}[htp]
    \centering
    \includegraphics[width=0.74\textwidth]{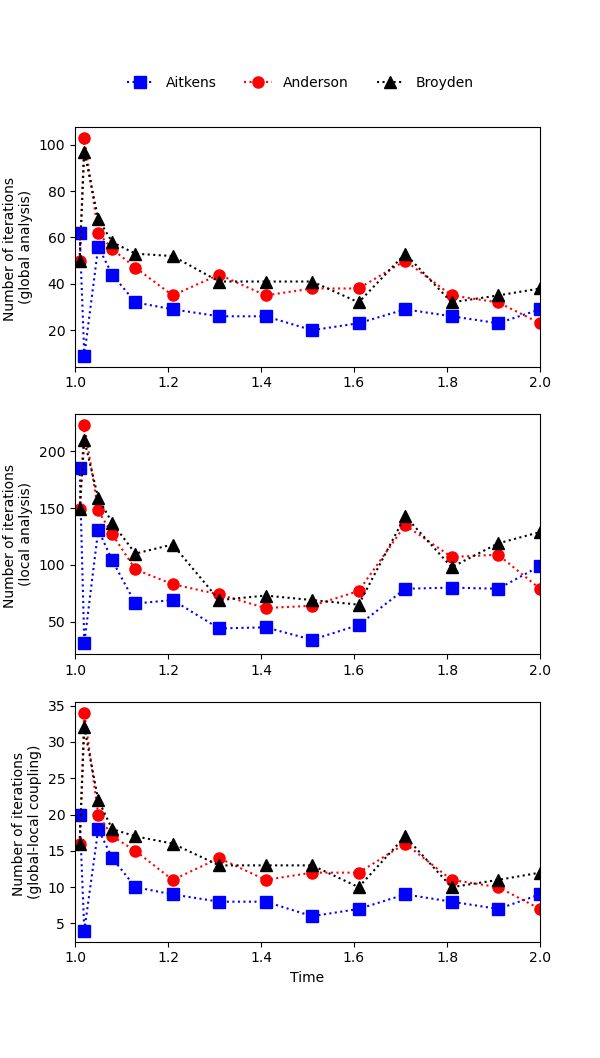}
    \caption{Comparison of the number of iterations in time between the GLIC-Aitken, GLIC-Anderson and GLIC-Broyden simulations for the bolted joint use case.}
    \label{fig:pwbiteracc}
\end{figure}

In terms of performance, the numbers of iterations in time are reported in Figure~\ref{fig:pwbiteracc} for the second step of the analysis only. The plot on the top shows the total number of Newton-Raphson iterations in the global analysis (elastic finite strain model with coarse uniform mesh), the plot in the middle shows the total number of Newton-Raphson iterations in the local analysis (where contact occurs) and the plot on the bottom shows the number of iterations for the global-local coupling. The blue squares, the red circles and the black triangles are results from GLIC-Aitken, GLIC-Anderson and GLIC-Broyden, respectively. This time, the red and black curves are not perfectly overlapped but are roughly on top of each other, whereas the blue curve is consistently on the bottom part of the plot. 

This means that the GLIC-Aitken turns out to be the faster accelerator technique when counting the number of total iterations during the analysis.
More precisely, as evaluated with numbers in Table~\ref{tab:perfo2}, global and local Newton-Raphson convergence rates are similar, but the global-local iteration is faster with GLIC-Aitken than with GLIC-Anderson and GLIC-Broyden. The number of increments is identical in the three cases, meaning that the convergence control heuristics work in the same way. The total number of Newton-Raphson iterations in the global analysis $N^G_\text{iter}$ is $1.49$ times smaller in GLIC-Aitken than in GLIC-Anderson and $1.59$ times than GLIC-Broyden. The total number of Newton-Raphson iterations in the local analysis $N^L_\text{iter}$ is $1.4$  times smaller in GLIC-Aitken than in GLIC-Anderson and $1.51$ times than in GLIC-Broyden.

\begin{table}[ht]
\centering
\begin{tabular}{lccccc} 
  & $N^G_\text{inc}$ & $N^G_\text{iter}$ & $N^L_\text{inc}$ & $N^L_\text{iter}$ & $N^{GL}$ \\
 \hline
 GLIC-Aitken   & $14$ & $434$ & $14$ & $1093$ & $151$ \\ 
 GLIC-Anderson & $14$ & $647$ & $14$ & $1533$ & $220$ \\
 GLIC-Broyden  & $14$ & $691$ & $14$ & $1648$ & $234$ \\
\end{tabular}
\caption{Comparison of total performance numbers between GLIC-Aitken, GLIC-Anderson and GLIC-Broyden for the bolted-joint use case.}
\label{tab:perfo2}
\end{table}

Relaxing the convergence criteria in the very same way of the first use case of Section~\ref{sec:case1} once more does not affect accuracy, as contact pressure and contact shear in $Z$ direction do not change, but improve performance.

\begin{figure}[htp]
    \centering
    \begin{subfigure}{0.475\textwidth}
        \includegraphics[width=\textwidth]{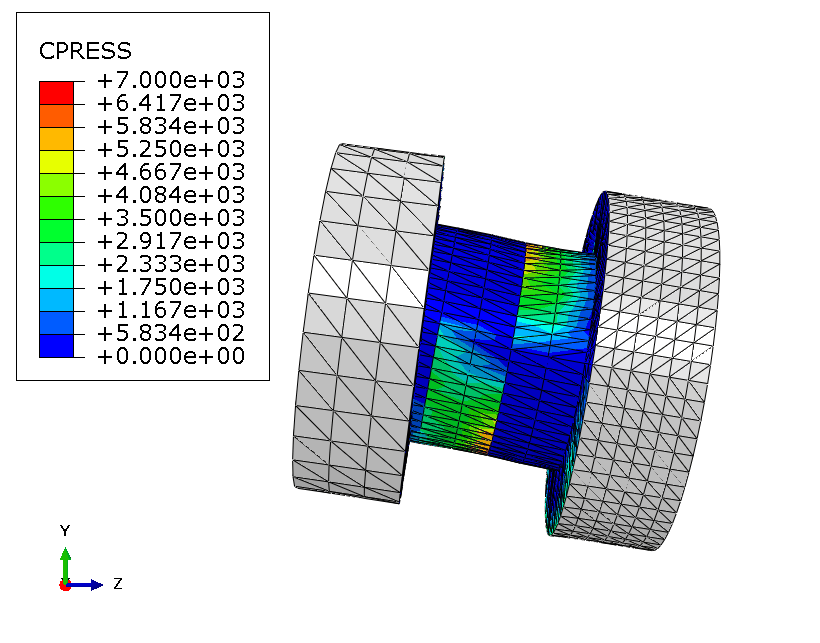}
        \caption{Contact pressure.}
    \end{subfigure}
    \begin{subfigure}{0.475\textwidth}
        \includegraphics[width=\textwidth]{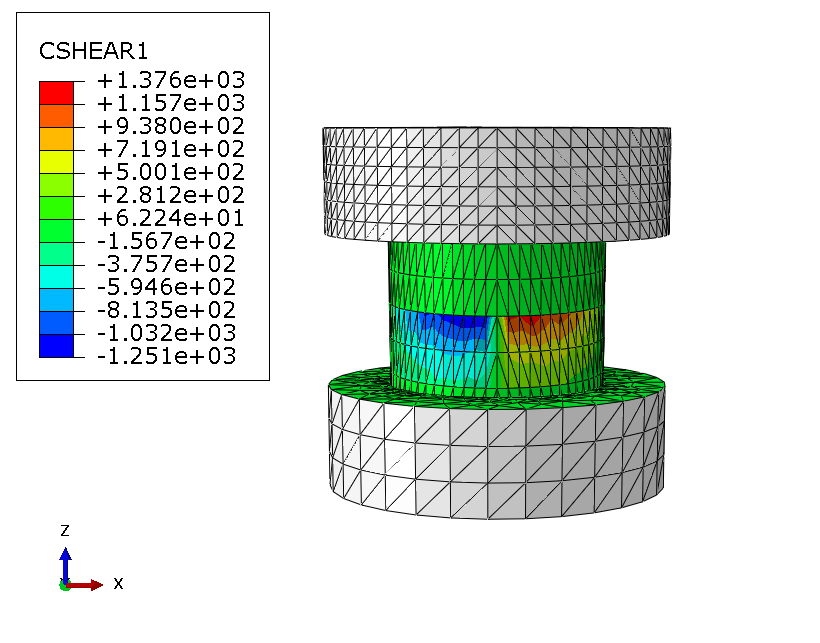}
        \caption{Contact shear ($Z$ direction).}
    \end{subfigure}
    \caption{Results in terms of contact interaction (MPa) with the GLIC-IS-Aitken simulation.}
    \label{fig:pwbcontactaitkensctrls}
\end{figure}

Figure~\ref{fig:pwbcontactaitkensctrls} shows contour plots of contact interaction results in terms of contact pressure and contact shear in $Z$ direction. Values are in line with reference results, maximum and minimum values being in between reference and co-simulation with Abaqus default convergence criteria.

\begin{figure}[htp]
    \centering
    \includegraphics[width=0.74\textwidth]{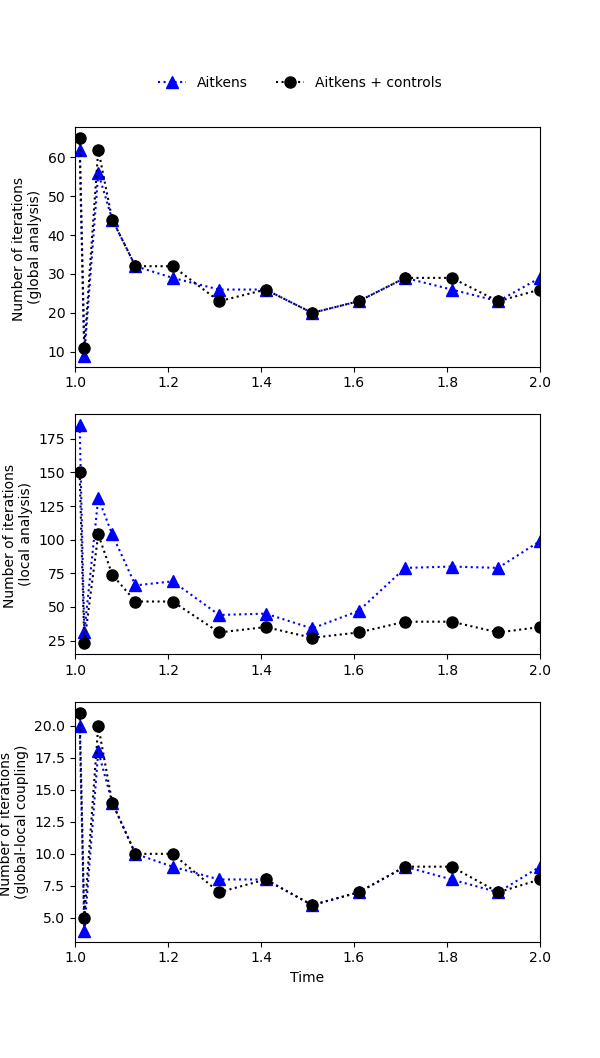}
    \caption{Comparison of the number of iterations in time between the GLIC-Aitken and GLIC-IS-Aitken simulations for the bolted joint use case.}
    \label{fig:pwsiterctrlsaitkens}
\end{figure}

In terms of performance, Figure~\ref{fig:pwsiterctrlsaitkens} shows the comparison in the total number of Newton-Raphson iterations in time between the Abaqus default convergence criteria and the relaxed convergence criteria. The performance gain using the latter strategy are appreciated chiefly on the local analysis, especially in the second half of the analysis, whereas there are no remarkable gains in the global analysis.

Also in this case, more precise numbers are listed in Table~\ref{tab:perfois2}, comparing GLIC-Aitken and GLIC-IS-Aitken. As curves of global Newton-Raphson iterations and global-local iterations were almost completely overlapped, the total number of global-local iterations $N^{GL}$ and the total number of Newton-Raphson iterations in the global analysis $N^G_\text{iter}$ are also in this case slightly smaller in GLIC-Aitken than in GLIC-IS-Aitken, namely, $2.6\%$ and $2.5\%$ smaller, respectively. On the other hand, the total number of Newton-Raphson iterations in the local analysis $N^L_\text{iter}$ is $1.5$ times smaller in GLIC-IS-Aitken than in GLIC-Aitken. Also in this case, the number of increments $N^G_\text{inc}$ and $N^L_\text{inc}$ are identical in GLIC-Aitken and GLIC-IS-Aitken.

\begin{table}[ht]
\centering
\begin{tabular}{lccccc} 
  & $N^G_\text{inc}$ & $N^G_\text{iter}$ & $N^L_\text{inc}$ & $N^L_\text{iter}$ & $N^{GL}$ \\
 \hline
 GLIC-Aitken    & $14$ & $434$ & $14$ & $1093$ & $151$ \\ 
 GLIC-IS-Aitken & $14$ & $445$ & $14$ & $727$ & $155$ \\
\end{tabular}
\caption{Comparison of total performance numbers between GLIC-Aitken and GLIC-IS-Aitken for the bolted joint use case.}
\label{tab:perfois2}
\end{table}

As also in this case the local model ($33393$ variables) is much larger than the global one ($8412$ variables), the GLIC-IS-Aitken approach is more convenient.

\section{Conclusions and future works}
\label{sec:conclusions}
Recent research works have proven the GLIC in conjunction with acceleration techniques as a robust and accurate solution for bridging multiple levels of modeling abstraction in computational structural mechanics, which is often a necessary approach when models are complex, in order to reduce the high and non-scalable modeling times and burdens. However, the simulation performance of the GLIC has always been reason of concern for the daily engineering practice. Within this article, the Authors study the GLIC implemented in Abaqus through co-simulation and compare acceleration techniques such as the Aitken's relaxation, Anderson acceleration and Broyden's quasi-Newton method. For the two highly nonlinear use cases under study, the Aitken's relaxation proved slightly more effective than the other techniques, as long as straightforward on implementation concerns The coupling reached acceptable levels of convergence with approximately $10$ iterations, in line with the most challenging multi-physics coupling strategies. Furthermore, the Authors could prove the employment of inexact solver strategy (GLIC-IS) more effective, with negligible loss of accuracy for engineering purposes, improving performance leveraging loosened tolerance values on the Newton-Raphson iteration schemes of the global and local analyses.

As the mitigation of the robustness and performance issues of the GLIC improves the appeal of the methodology for an effective engineering practice, the full promise of non-intrusive modeling is not fully met yet: in the presented works, the authors needed to manufacture local models with solid-to-shell coupling in order to align global and local interfaces to the facets of global elements, in order to avoid to cut through global element with the local patch. This issue was successfully tackled with a re-meshing strategy in the local model by \cite{blanchard_2019}, translating modeling intrusivity to implementation intrusivity. The authors would also like to explore new opportunities in future works such as the employment of virtual elements \cite{wriggers_2023} or enhanced elements \cite{hauseux:hal-01962742} where global elements are cut by the coupling interface.

\section*{Aknowledgements}
The Authors thanks Albert Kurkchubasche and Sethu Subbarayalu for their precious help on the implementation tasks within the Abaqus co-simula\-tion tools.





\bibliographystyle{plain} 
\bibliography{papers}





\end{document}